\documentclass[11pt]{amsart}
\usepackage{amssymb, amstext, amscd, amsmath}
\makeatletter
\def\@cite#1#2{{\m@th\upshape\bfseries%
[{#1\if@tempswa{\m@th\upshape\mdseries, #2}\fi}]}}
\makeatother

\theoremstyle{plain}
\newtheorem{thm}{Theorem}[section]
\newtheorem{cor}[thm]{Corollary}
\newtheorem{prop}[thm]{Proposition}
\newtheorem{lem}[thm]{Lemma}

\theoremstyle{definition}
\newtheorem{rem}[thm]{Remark}

\newtheorem{defn}[thm]{Definition}

\newtheorem{eg}[thm]{Example}
\newtheorem{ques}[thm]{Question}

\newcommand{\bC}{{\mathbb{C}}}

\newcommand{\bN}{{\mathbb{N}}}
\newcommand{\bQ}{{\mathbb{Q}}}
\newcommand{\bR}{{\mathbb{R}}}

  \newcommand{\B}{{\mathcal{B}}}
  \newcommand{\C}{{\mathcal{C}}}
  \newcommand{\D}{{\mathcal{D}}}

\renewcommand{\H}{{\mathcal{H}}}
  \newcommand{\I}{{\mathcal{I}}}

\renewcommand{\L}{{\mathcal{L}}}
  \newcommand{\M}{{\mathcal{M}}}
  \newcommand{\N}{{\mathcal{N}}}
\renewcommand{\O}{{\mathcal{O}}}
\renewcommand{\P}{{\mathcal{P}}}
  
  \newcommand{\R}{{\mathcal{R}}}

  \newcommand{\T}{{\mathcal{T}}}
  \newcommand{\U}{{\mathcal{U}}}

\newcommand{\ep}{\varepsilon}
\renewcommand{\phi}{\varphi}
\newcommand{\upchi}{{\raise.35ex\hbox{$\chi$}}}

\newcommand{\fA}{{\mathfrak{A}}}
\newcommand{\fB}{{\mathfrak{B}}}

\newcommand{\fD}{{\mathfrak{D}}}

\newcommand{\fJ}{{\mathfrak{J}}}

\newcommand{\fM}{{\mathfrak{M}}}

\newcommand{\fX}{{\mathfrak{X}}}

\newcommand{\AND}{\text{ and }}
\newcommand{\FOR}{\text{ for }}
\newcommand{\FORAL}{\text{ for all }}

\newcommand{\qand}{\quad\text{and}\quad}

\newcommand{\qfor}{\quad\text{for}\quad}
\newcommand{\qforal}{\quad\text{for all}\quad}

\newcommand{\diag}{\operatorname{diag}}

\newcommand{\id}{\operatorname{id}}

\newcommand{\rank}{\operatorname{rank}}

\newcommand{\spn}{\operatorname{span}}

\newcommand{\norm}[1]{\left\| #1 \right\|}
\newcommand{\ol}{\overline}
\newenvironment{sbmatrix}{\left[\begin{smallmatrix}}{\end{smallmatrix}\right]}

\newcommand{\wot}{\textsc{wot}}
\DeclareMathOperator*{\wotlim}{\textsc{wot}--lim}

\newcommand{\ltsr}{\mathrm{ltsr}}
\newcommand{\rtsr}{\mathrm{rtsr}}
\newcommand{\tsr}{\mathrm{tsr}}

\newcommand{\AD}{\mathrm{A}(\mathbb{D})}

\newcommand{\linf}{\ell^\infty }
\newcommand{\ltwo}{\ell^2}

\title[Topological stable rank]{Topological stable rank of nest algebras}

\author[K.R.Davidson]{Kenneth R. Davidson}
\address{Pure Math.\ Dept.\\U. Waterloo\\Waterloo,
ON\; N2L--3G1\\CANADA} \email{krdavids@uwaterloo.ca}
\thanks{First author partially supported by an NSERC grant.}

\author[Y.Q. Ji]{You Qing Ji}
\address{Math. Dept. \\ Jilin University\\ Changchun 130012\\P.R. CHINA}
\email{jiyq@jlu.edu.cn}
\thanks{Second author partially supported by RFDP(20050183002), 
NCET and the China Scholarship Council.}


\date{}
\begin{document}

\maketitle

\begin{abstract}
We establish a general result about extending a right invertible 
row over a Banach algebra to an invertible matrix.
This is applied to the computation of right topological stable rank
of a split exact sequence.
We also introduce a quantitative measure of stable rank.
These results are applied to compute the right (left) topological 
stable rank for all nest algebras.  
This value is either 2 or infinity, and $\rtsr(\T(\N)) = 2$ occurs only
when $\N$ is of ordinal type less than $\omega^2$ and the dimensions
of the atoms grows sufficiently quickly.
We introduce general results on `partial matrix algebras'
over a Banach algebra.  This is used to obtain an inequality
akin to Rieffel's formula for matrix algebras over a Banach algebra.
This is used to give further insight into the nest case.
\end{abstract}

\section{Introduction}

The topological stable rank for Banach algebras was defined in Rieffel~\cite{R} 
as a non-commutative analogue of the covering dimension for compact spaces
that was modelled on the Bass stable rank of rings \cite{B}.
In \cite{DLMR}, the first example of an operator algebra for which the left and right
topological stable ranks differ was produced.  Indeed a large class of nest algebras
were shown to have this property.
The analysis of topological stable rank for arbitrary nests was begun there.  
The purpose of this paper is to develop new techniques for computing the topological
stable rank of Banach algebras.  The particular application we focus on is the
calculation of the topological stable rank of all nest algebras.

Given a unital Banach algebra $\fA$, we denote by $Rg_n(\fA)$ (resp.\ $Lg_n(\fA))$ 
the set of $n$-tuples of elements of $\fA$ which generate $\fA$ as a right ideal 
(resp.\ as a left ideal).   
That is,
\[
 Rg_n(\fA) = \{ (a_1, \dots, a_n) \in \fA^n : 
 \exists (b_1, \dots, b_n) \in \fA^n 
 \text{ with } \sum_{i=1}^n a_i b_i = 1 \} . 
\]
The \textit{right topological stable rank} of $\fA$, 
denoted by $\rtsr(\fA)$, is the least positive integer $n$
for which $Rg_n(\fA)$ is dense in $\fA^n$.   
When no such integer exists, we set $\rtsr(\fA) = \infty$. 
The \textit{left topological stable rank} of $\fA$, $\ltsr(\fA)$, is defined analogously.
If $\ltsr(\fA) = \rtsr(\fA)$, we refer to their common value simply as the 
\textit{topological stable rank} of $\fA$, written $\tsr(\fA)$.   
When $\fA$ is not unital, we define the right (left) topological stable rank of 
$\fA$ to be that of its unitization.

When $\fA$ is an operator algebra acting on a Hilbert space $\H$,
it is convenient to reformulate these notions in terms of the 
row and column spaces of $\fA$.  
Let $\R_n(\fA)$ denote the $1 \times n$ matrices 
$A = \big[ A_1\ \dots\ A_n \big]$ with coefficients in $\fA$,
considered as a subspace of $\B(\H^{(n)},\H)$ with the induced norm.
Likewise let $\C_n(\fA)$ denote the space of $n \times 1$ matrices
with coefficients in $\fA$ normed as a subspace of $\B(\H,\H^{(n)})$.
Then we identify $Rg_n(\fA)$ with the right invertible elements of $\R_n(\fA)$,
and $Lg_n(\fA)$ with the left invertible elements of $\C_n(\fA)$.
Density can now be considered with respect to this operator space norm.

When manipulating elements of $Rg_n(\fA)$, it is useful to observe that
multiplication on the right by an invertible $n \times n$ matrix preserves this set.
So it is helpful to be able to produce an ample supply of such matrices.
In particular, it is especially convenient if there is an invertible matrix
with prescribed first row and prescribed first column of the inverse.
Being able to do this is a strong finiteness condition, and we call
an algebra with this property \textit{completely finite}.
We characterize when this is possible in general Banach algebras, 
and determine exactly which nest algebras have this property.

Rieffel establishes several results about ideals and quotients.
However his hypotheses require the ideal to have a bounded
approximate identity.  In practice, this rarely occurs in the study
of nonself-adjoint operator algebras such as nest algebras.
We have found that a suitable replacement is to assume that the quotient map splits.
This occurs, for example, with nest algebras when the ideal is the
kernel of an expectation onto the block diagonal algebra corresponding
to a subnest.  The added bonus is that if the quotient is also 
completely finite, then we obtain a precise formula for the stable rank
of the algebra:
\[ \rtsr(\fA) = \max\{ \rtsr(\fA/\fJ), \rtsr(\fJ) \} .\]

Then we introduce a quantitative measure of topological stable rank.
This is crucial for comparing the topological stable rank  of an
$\linf$ direct sum of Banach algebras.  These algebras occur frequently
as the quotients of nest algebras as mentioned in the previous paragraph.
In particular, it is necessary to get a quantitative measure in order to
establish lower bounds for the topological stable rank.
One method used in \cite{DLMR} to show, for example, that the usual
algebra of upper triangular operators $\T(\bN)$ has $\rtsr(\T(\bN)) = \infty$
was to use a result of \cite{DHO} that exhibits a surjective homomorphism
of $\T(\bN)$ onto $\T(\bN)^*$.  Unfortunately this method is limited in that
it requires some uniform behaviour of the nest (namely that it has arbitrarily long
disjoint intervals which have all rank 1 atoms).  This is too stringent.
The authors of \cite{DLMR} tried hard to find this quantitative obstruction
but failed.  Theorem~\ref{T:lower bound} and its corollaries set out a
general method for establishing that $\rtsr(\T(\N)) = \infty$ whenever
the dimensions of the atoms do not grow too fast.

In \cite{DLMR}, it was shown that many nests have 
\[ \rtsr(\T(\N)) = \ltsr(\T(\N)) = \infty .\]
This includes all uncountable nests and any nest with an infinite rank atom.
Indeed, the problem was reduced to computing $\rtsr(\T(\N))$ when $\N$
is of ordinal type with finite dimensional atoms.
We show here that if the ordinal is at least $\omega^2$, then again
$\rtsr(\T(\N)) = \ltsr(\T(\N)) = \infty$.
Since a finite upper triangular matrix of Banach algebras has topological stable
rank determined by the diagonal entries, the problem is therefore reduced to
studying nests of order type $\omega$ with finite rank atoms.

It is evident for C*-algebras that the left and right topological stable ranks agree.
Herman and Vaserstein \cite{HV} showed that for C*-algebras, the topological
stable rank coincides with the Bass stable rank.
However this is not true for all Banach algebras.  For example \cite{JMW}, 
the disk algebra has $\tsr(\AD) = 2$ but the Bass stable rank is 1.
Vaserstein \cite{V} (also Warfield \cite{W}) showed that the left and right
Bass stable rank coincide for all rings.
So it was a surprising result from \cite{DLMR}  that the 
left and right topological stable ranks of an operator algebra can differ:

\begin{thm}\cite[Theorem 2.11]{DLMR}  \label{DLMR2.11}
Let $\N$ be a nest of order type $\omega$ with finite rank atoms
of dimension $n_k$ for $k\ge1$. 
Suppose that there exists an $r > 1$ and an integer $J \ge1$  such that 
$\displaystyle \max_{1 \le i \le  (k+1)J} n_i \ge r  \max_{1 \le i \le  kJ} n_i$  for all $k \ge 1$. 
Then  $\ltsr(\T(\N)) = \infty$ and  $\rtsr(\T(\N)) = 2$.
\end{thm}

When the ranks of the atoms is monotone (or even monotone within a constant),
we obtain a very simple way to compute the topological stable rank.
Let $d_j$ count the number of atoms with dimension in the interval
$(2^{j-1},2^j]$.  If $\sup d_j = \infty$, then $\rtsr(\T(\N)) = \infty$.
While if  $\sup d_j < \infty$, then $\rtsr(\T(\N)) = 2$.
This is a significant sharpening of the results of \cite{DLMR}.
Indeed, we provide a simpler proof of the theorem above,
by reducing it to the simpler argument of \cite[Theorem~2.1]{DLMR}
which assumes rapid exponential growth $n_k \ge 4 \sum_{i<n} n_i$.
Indeed, our methods produce an algorithm for computing the topological stable
rank when the nest has a finite sequence of subnests which each have a
monotonicity property relative to the next nest.

While this sounds like a rather special case, we show that it applies whenever
the right topological stable rank is finite.  We obtain a quantitative invariant
which is equivalent to \mbox{$\rtsr(\T(\N)) < \infty$.}  Combining this with the
results on relative monotonicity, we are able to conclude that $\rtsr(\T(\N)) = 2$
in this case.  So it follows that for all nest algebras, $\rtsr(\T(\N))$ is either $2$ or $\infty$.
The computation of this invariant may be difficult.  But the proof establishes
a method for computing $\rtsr(\T(\N))$ based on the ideas in the previous section.

Rieffel established a remarkable precise formula for the topological stable
rank of a matrix algebra over $\fA$, namely
\[ \rtsr(\fM_n(\fA)) = \left\lceil \frac{\rtsr(\fA) - 1}n \right\rceil + 1 .\]
In particular, for $n$ sufficiently large, $\rtsr(\fM_n(\fA))$ takes only the
values 1,2 and $\infty$.
The case of $\rtsr(\fA)=1$ was shown by Rieffel to be equivalent to the
density of the invertible elements in $\fA$; and so implies that $\ltsr(\fA)=1$ also.
This never happens in nest algebras because either $\T(\N)$ or $\T(\N)^*$
has a quotient which has two isometries with orthogonal ranges;
and hence $\rtsr(\T(\N)) = \infty$ or $\ltsr(\T(\N)) = \infty$.
In the case of a nest of order type $\omega$, it is the latter which occurs.
In the last section, we consider `partial matrix algebras' over a Banach algebra,
and obtain an inequality analogous to Rieffel's formula.
This allows us to show that the stable ranks
of two nest algebras are equal whenever the dimensions of their atoms are related 
by an inequality $cn_k \le m_k \le d m_k$ for $0 < c \le d < \infty$.

We briefly remind the reader of the basic notation for nests.
A nest $\N$ is a complete chain (with respect to inclusion) of closed subspaces 
of a Hilbert space $\H$ containing $\{0\}$ and $\H$.
The nest algebra $\T(\N)$ is the algebra of all operators $T\in\B(\H)$
which leave each element $N \in \N$ invariant.
An interval of $\N$ is a subspace $N\ominus M$ for elements $M<N \in \N$.
An atom is a minimal interval. 
For $N \in \N$, let $N^- = \bigvee\{N' \in \N : N' < N\}$.
We write 
\[ \T_0(\N) = \{T \in \T(\N) : P_ATP_A=0 \FORAL \text{atoms }A \} \] 
for the ideal of strictly upper triangular operators.

When $\N$ is well ordered, we say that $\N$ has order type $\alpha$ if the 
atoms of $\N$ are of order type $\alpha$.  So a nest of order type $\omega$
has the form $\N=\{ N_k, \H : k \ge 0\}$ with $N_{k-1} < N_k$
and  atoms $A_k = N_k \ominus N_{k-1}$ for $k \in \bN$.
See \cite{D_nest} for further background.

\section{Invertible Matrices} \label{S:invertible}

It will be useful to have a criterion that allows us to complete
the first row of a matrix to an invertible one.
The natural condition on a row $R = \big[ A_1\ \dots \ A_n \big] \in \R_n(\fA)$
is the existence of a right inverse $C = \big[ C_1\ \dots \ C_n \big]^t \in \C_n(\fA)$
such that $RC=1$.
When $n=1$, this requires that whenever $A$ is right invertible, then
it is invertible---a finiteness condition on $\fA$.

Likewise, one can see that $\B(\H)$ does not have this property for any $n$
because of its infinite nature.
Let $S_1,\dots,S_n$ be $n$ isometries with pairwise orthogonal ranges
summing to the whole space (generators of $\O_n$).
Then $S = \big[ S_1\ \dots \ S_n \big]$ has right inverse $S^*$.
However, as $S$ is an isometry, it is injective and so there
are no non-zero columns $C$ satisfying $SC=0$.

On the other hand, this is an easy exercise for the operator algebra $\fM_k$.
For then the condition that $RC= I_k$ guarantees that the rank of $R$
is exactly $k$.  Hence one can easily extend $R$ to a $(nk)\times (nk)$ matrix
of full rank, whence invertible in $\fM_n(\fM_k)$.

We first establish a general condition, and then specialize to
those nest algebras which are finite in this sense.
Let $E_{ij}$ denote the standard matrix units of $\fM_n(\fA)$.

\begin{thm} \label{T:extend}
Let $\fA$ be a unital Banach algebra, and let $n \ge 2$.  Suppose that 
$R = \big[ A_1\ \dots \ A_n \big] \in \R_n(\fA)$
and $C = \big[ C_1\ \dots \ C_n \big]^t \in \C_n(\fA)$
satisfy $RC=1$. Then the following are equivalent:
\begin{enumerate}
\item There is an invertible matrix $W \in \fM_n(\fA)$ with first row $R$.
\item There is an invertible matrix $W \in \fM_n(\fA)$ with first column $C$.
\item There is an invertible matrix $W \in \fM_n(\fA)$ with first row $R$
 such that the first column of $W^{-1}$ is $C$.
\item The idempotent $CR$ is similar in $\fM_n(\fA)$ to $E_{11}$.
\end{enumerate}
\end{thm}

\begin{proof}
Evidently (iii) implies (i) and (ii).
Conversely, suppose that (i) holds with invertible matrix
$W = \begin{bmatrix}\ R\ \\ X \end{bmatrix}$, where $X \in \fM_{n-1,n}(\fA)$.
Write $W^{-1} = \begin{bmatrix}D & Y \end{bmatrix}$, 
where $D \in \C_n(\fA)$ and $Y \in \fM_{n,n-1}(\fA)$.
Then $RY=0$ and $XY = I_{n-1}$.
Let $V =  \begin{bmatrix}C & Y \end{bmatrix}$.  Observe that
\[
 WV = \begin{bmatrix} 1 & 0\\ XC & I_{n-1} \end{bmatrix} 
\]
is invertible.  Hence 
\[
 V^{-1} = (WV)^{-1}W =
 \begin{bmatrix} 1 & 0  \\ -XC & I_{n-1} \end{bmatrix}
 \begin{bmatrix} \ R\  \\ X \end{bmatrix}  =
 \begin{bmatrix} \ R \ \\ X' \end{bmatrix}
\]
Thus (iii) holds.  The proof that (ii) implies (iii) is the same.\vspace{.3ex}

Now assume (iii), and let $W  = \begin{bmatrix}\ R\ \\ X \end{bmatrix}$
with inverse $W^{-1} = \begin{bmatrix} C & Y \end{bmatrix}$.
Then 
\[
 RW^{-1} = \big[ 1 \ 0 \ \dots \ 0 \big] \qand  
 WC = \big[ 1 \ 0 \ \dots \ 0 \big]^t .
\]
So $W(CR)W^{-1} = E_{11}$.

Conversely, suppose that (iv) holds and $W(CR)W^{-1} = E_{11}$.
Write $W  = \begin{bmatrix}\ R'\ \\ X \end{bmatrix}$
and $W^{-1} = \begin{bmatrix} C' & Y \end{bmatrix}$,
where $R'\in\R_n(\fA)$ and $C'\in\C_n(\fA)$.
Then $CR = W^{-1}E_{11}W = C'R'$ and $R'C'=1$.
Let $A=RC'$ and $B=R'C$; and note that 
\[ AB = R(C'R')C=RCRC = 1 \]
and 
\[ BA = R'(CR)C' = R'C'R'C'=1 . \]
Moreover, 
\[ AR' = R(C'R') = RCR = R \]
and 
\[ C'B = (C'R')C = CRC = C .\]
Therefore
\[ 
 \begin{bmatrix} A&0\\0&I_{n-1} \end{bmatrix}
 \begin{bmatrix}\ R' \\  X \end{bmatrix} =
 \begin{bmatrix}\ R \ \\ X \end{bmatrix}
\]
is an invertible matrix with first row $R$; and its inverse is
\[ 
 \begin{bmatrix}C'&Y \end{bmatrix} 
 \begin{bmatrix} B&0\\0&I_{n-1} \end{bmatrix} =
 \begin{bmatrix}C&Y \end{bmatrix}
\]
So (iii) holds.
\end{proof}

\begin{defn}
A unital Banach algebra $\fA$ is \textit{completely finite}
if every right invertible $R \in \R_n(\fA)$ is the first row of an invertible 
matrix in $\fM_n(\fA)$ for all $n \ge1$.

A Banach algebra $\fA$ is \textit{finite} if every right (left) invertible
element is invertible; 
and $\fA$ is \textit{stably finite} if $\fM_n(\fA)$
is finite for all $n \ge 1$.
\end{defn}

\begin{thm} \label{T:complete rectangles}
Suppose that $\fA$ is completely finite.
If $1 \le k \le n$ and $A = \big[ a_{ij} \big] \in \fM_{k \times n}(\fA)$
has a right inverse $B \in \fM_{n \times k}$, then there is
an invertible matrix $W \in \fM_n(\fA)$ with the first $k$ rows equal to $A$
and such that the first $k$ columns of $W^{-1}$ equal $B$.
\end{thm}

\begin{proof}
We proceed by induction on $k$.  
By hypothesis, the result is valid for $k=1$.
Assume that $2 \le k \le n$ and that we have established the result for $k-1$.

Write $A = \begin{bmatrix}A_1\\ R \end{bmatrix}$ where $A_1$
is the first $k-1$ rows of $A$ and $R$ is the $k$th row.
Similarly write $B =  \begin{bmatrix}B_1 & C \end{bmatrix}$.
We are given that $AB = I_k$.
Then since $A_1B_1 = I_{k-1}$, the assumption is that there
is an invertible matrix $W_1 = \begin{bmatrix}A_1\\ A_2 \end{bmatrix}$
in $\fM_n(\fA)$ with inverse of the form 
$W_1^{-1} = \begin{bmatrix}B_1 & B_2 \end{bmatrix}$.
Then 
\[
  AW_1^{-1} = \begin{bmatrix}A_1\\ R \end{bmatrix}
  \begin{bmatrix}B_1 & B_2 \end{bmatrix} =
  \begin{bmatrix}I_{k-1} & 0\\ * & R' \end{bmatrix}
\]
and
\[
  W_1B = \begin{bmatrix}A_1\\ A_2 \end{bmatrix}
  \begin{bmatrix}B_1 & C \end{bmatrix} =
  \begin{bmatrix}I_{k-1} & * \\ 0 &C' \end{bmatrix} .
\]
Since $I_k = (AW_1^{-1})(W_1B)$, it follows that both of the entries 
marked $*$ are $0$.

Now we have $R'C' = I$ where 
$R' \in \R_{n+1-k}(\fA)$ and $C' \in \C_{n+1-k}(\fA)$.
As $\fA$ is completely finite, there is an invertible matrix 
$X \in \fM_{n+1-k}(\fA)$ with first row $R'$ and so that $X^{-1}$
has first column $C'$.
Therefore $W_2 = I_{k-1}\oplus X$ is invertible in $\fM_n(\fA)$.
We have
\begin{align*}
 AW_1^{-1}W_2^{-1} &= 
 \begin{bmatrix}I_{k-1} & 0\\ 0 & R' \end{bmatrix}
 \begin{bmatrix}I_{k-1} & 0\\ 0 & X^{-1} \end{bmatrix} \\ &=
 \begin{bmatrix}I_{k-1} & 0\\ 0 &  \begin{sbmatrix}I& 0\end{sbmatrix} \end{bmatrix} = 
 \begin{bmatrix}I_k & 0 \end{bmatrix} 
\end{align*}
and
\begin{align*}
  W_2W_1B &= \begin{bmatrix}I_{k-1}&0\\ 0& X \end{bmatrix}
  \begin{bmatrix}I_{k-1} & 0 \\ 0 & C' \end{bmatrix} \\ &=
  \begin{bmatrix}I_{k-1} & 0 \\ 0 & \begin{sbmatrix}I\\ 0\end{sbmatrix} \end{bmatrix} =
  \begin{bmatrix}\ I_k\  \\ 0 \end{bmatrix} .
\end{align*}

Let $W = W_2W_1$. 
Observe that if we write $W = \begin{bmatrix}Y\\ Z \end{bmatrix}$ where
$Y$ consists of the first $k$ rows, then
\[ 
 A = AW^{-1}W = 
 \begin{bmatrix}I_k & 0 \end{bmatrix}  
 \begin{bmatrix} Y\\ Z \end{bmatrix}  
 = Y .
\]
So the first $k$ rows of $W$ are precisely $A$. 
Similarly the first $k$ columns of $W^{-1}$ equals $B$.
\end{proof}

\begin{cor} \label{comp finite is stably finite}
If $\fA$ is completely finite, then it is stably finite.
\end{cor}

\begin{proof}
Let $A \in \fM_n(\fA)$ be right invertible.
In the case $k=n$ of Theorem~\ref{T:complete rectangles}, 
the invertible matrix agreeing with $A$
on the first $n$ rows must be $A$ itself.  
So $A$ is invertible.
\end{proof}

\begin{cor} \label{matrix comp finite}
If $\fA$ is completely finite, then $\fM_k(\fA)$ is completely finite for all $k \ge 1$.
\end{cor}

\begin{proof}
A right invertible row in $\R_n(\fM_k(\fA))$ is a right invertible 
matrix in \mbox{$\fM_{k\times kn}(\fA)$.}  So it may be completed to
an invertible matrix in $\fM_{kn}(\fA)$ by Theorem~\ref{T:complete rectangles}.
\end{proof}

We make the following elementary observation for future use:

\begin{prop}\label{comp_fin_quot}
Suppose that $\fA$ is completely finite.
If $\fJ$ is an ideal of $\fA$ such that the quotient map splits,
then  $\fA/\fJ$ is completely finite.
\end{prop}

\begin{proof}
Suppose that $\pi$ is the quotient map and $\sigma:\fA/\fJ \to \fA$ is the splitting map.
If $R = \big[ A_1\ \dots \ A_n \big] \in \R_n(\fA/\fJ)$ has a  
right inverse $C \in \C_n(\fA/\fJ)$, then 
$R' = \sigma(R) =  \big[ \sigma(A_1)\ \dots \ \sigma(A_n) \big]$ in $\R_n(\fA)$
has right inverse $\sigma(C)$.
Use the fact that $\fA$ is completely finite to find an invertible operator 
$W \in \fM_n(\fA)$ with first row $R'$, say $W = \big[ W_{ij} \big]$.
Then it is easily checked that $\pi(W) = \big[ \pi(W_{ij}) \big]$ is the
desired invertible element of $\fM_n(\fA/\fJ)$ with first row equal to $R$.
\end{proof}

We now specialize to nest algebras.
If a nest $\N$ contains an atom of infinite rank, then there is a homomorphism
of $\T(\N)$ onto $\B(\H)$.  The comments preceding the proof show that 
a right invertible row cannot, in general, be extended to an invertible matrix.
Also, if $\N$ is a continuous nest acting on a separable Hilbert space, then
by the Similarity Theorem \cite{D_sim}, then $\N$ is similar to its infinite ampliation;
so $\T(\N) \sim \T(\N)\otimes \B(\H')$.
Then using $n$ isometries $S_1,\dots,S_n$ with orthogonal ranges in $\H'$
summing to the whole space, we can define the row operator
\[ R = \big[ I \otimes S_1\ I\otimes S_2\ \dots \ I\otimes S_n \big] .\]
Then $R$ is an isometry in $\R_n(\T(\N)\otimes \B(\H'))$ such that 
$R^*$ belongs to $\C_n(\T(\N)\otimes \B(\H'))$.
So $RR^* = I$, but this does not extend to an invertible matrix.

Indeed, if $\N$ is uncountable, then it is similar to a nest with continuous part.
Thus one can see that it contains right invertible but non-invertible elements.
So $\T(\N)$ is not finite.  If $\N$ has an infinite rank atom $E$, or a non-atomic interval $E$,
then there is a quotient of $\T(\N)$ onto $\T(E \cap \N)$; and this map evidently splits.
So the proof of Proposition~\ref{comp_fin_quot} and the remarks of the previous paragraph
show that there are right invertible row operators of arbitrary size that do not extend
to invertible matrices.
However, as we will see below, even when $\T(\N)$ is not finite, it may
have this matrix extension property.

Theorem~\ref{T:extend} points to the issue of when two idempotents 
in a nest algebra are similar.  This problem has been completely solved
by Larson and Pitts \cite{LP}.  In this first result, we only need the result for
countable nests, where more elementary methods suffice.

\begin{thm} \label{T:extend_nest}
Let $\N$ be a nest on a separable Hilbert space.
Then $\T(\N)$ is completely finite if and only if $\N$ is
countable nest which has only finite rank atoms.
\end{thm}

\begin{proof}
The discussion prior to the theorem shows that the conditions
on $\N$ are necessary for $\T(\N)$ to be finite; so assume this form.
Let the atoms of $\N$ be enumerated as $E_i$ for $i \ge 1$,
and let their ranks be $n_i \in \bN$. 
Also let $N_i\in\N$ such that $E_i = N_i \ominus N_i^-$.

First let us deal with the case $n=1$.
Suppose that $A,B \in \T(\N)$ such that $AB=I$; and let $P=BA$.
Then it suffices to show that $P=I$, or equivalently that $Q=I-P = 0$.
Now $P=P^2$ and $Q=Q^2$ are idempotent.
The diagonal expectation $\Delta(X) = \sum_i P_{E_i}XP_{E_i}$ satisfies
\[ I = \Delta(I) = \Delta(A) \Delta(B) \qand \Delta(P) = \Delta(A) \Delta(B) .\]
Since the atoms are finite dimensional, the diagonal $\fD$ of $\T(\N)$
is a finite von Neumann algebra.  
It follows that $\Delta(B) = \Delta(A)^{-1}$ and $\Delta(P) = I$.
Therefore $\Delta(Q) = 0$.
In particular,  $QN_i \subset N_i$ and $P_{E_i}QP_{E_i}=0$ 
imply that $QN_i \subset N_i^-.$

There is a standard argument from the similarity theory of nests that 
now shows that $Q=0$.  Indeed, there is an invertible operator $S \in \B(\H)$
so that $Q' = SQS^{-1}$ is a self-adjoint projection.  
This belongs to the nest algebra $\T(S\N)$; and hence to the diagonal algebra.
Because $\N$ is countable, so is $S\N$.  Hence it has no non-atomic part.
So the expectation $\Delta'$ onto the diagonal of $\T(S\N)$ is given by
$\Delta'(X) = \sum_i P_{E'_i}XP_{E'_i}$, where
$E'_i = SN_i \ominus SN_i^-$ in $S \N$.
Now compute
\begin{align*}
 P_{E'_i} Q' P_{E'_i}\H &= P_{SN_i^-}^\perp SQS^{-1} P_{SN_i} \H\\
 &\subset P_{SN_i^-}^\perp SQ S^{-1}SN_i = P_{N_i^-}^\perp SQ N_i\\
 &\subset P_{SN_i^-}^\perp S N_i^- = \{0\} .
\end{align*}
Hence $Q' = \Delta'(Q') = 0$.  Consequently $Q=0$ and $\T(\N)$ is finite.

Now take $n \ge 2$.
Fix $R \in \R_n(\T(\N))$ and $C \in \C_n(\T(\N))$ such that $RC=1$.
The idempotent $P = CR$ in $\fM_n(\T(\N))$ may be considered as
an element of $\T(\N^{(n)})$; and likewise we consider $Q = E_{11}$ in this way.  
Since $P_{E_i}RCP_{E_i} = P_{E_i}$ has rank $n_i$, it follows that
the rank of $P_{E_i}^{(n)} CR P_{E_i}^{(n)} $ is also $n_i$.
This agrees with the rank of $P_{E_i}^{(n)} Q P_{E_i}^{(n)}$.

The nest on $P\H^{(n)} \oplus (I_n - P)\H^{(n)}$ consisting of the subspaces
$PN^{(n)} \oplus (I_n - P) N^{(n)}$ for $N \in \N$ is a countable nest which is
order isomorphic to $\N^{(n)}$ and this isomorphism preserves dimensions of atoms.
We consider $N^{(n)} \simeq N \oplus N^{(n-1)}$ 
where the first summand is identified with $QN^{(n)}$.
Because these nests are countable, they are unitarily equivalent via a unitary
which carries $PN^{(n)}$ onto $UPN^{(n)} = QN^{(n)}$ and carries 
$(I_n - P) N^{(n)}$ onto $U(I_n - P) N^{(n)} = Q^\perp N^{(n)}$.
Let $V$ be the invertible operator which takes $\H^{(n)}$ onto 
$P\H^{(n)} \oplus (I_n - P)\H^{(n)}$ by $Vx = Px \oplus (I_n - P)x$.
Then it is evident that $(UV) P (UV)^{-1} = Q$.
The result now follows from Theorem~\ref{T:extend}.
\end{proof}

Now we return to the issue of extending right invertible rows to invertible
matrices in other nests.  There is a somewhat surprising answer.
The class of nests in this next result includes many uncountable
nests, such as the Cantor nest on $\ltwo(\bQ)$ with basis $\{e_q : q \in \bQ\}$ 
consisting of the subspaces $N_r = \spn\{ e_q : q \le r \}$ for $r \in \bR$
and $N_r^- = \spn\{ e_q : q < r \}$ for $r \in \bQ$ together with $\{0\}$ and $\H$.
This has a dense set of one-dimensional atoms.
Since an uncountable nest $\N$ is similar to a nest with continuous part,
$\T(\N)$ is not finite.

\begin{thm} \label{T:extend_more_nests}
Let $\N$ be a nest on a separable Hilbert space. 
Then the following are equivalent:
\begin{enumerate}
\item For every $n\ge2$ and every right invertible 
 $R = \big[ A_1\ \dots \ A_n \big]$ in $\R_n(\T(\N))$
 with right inverse $C = \big[ C_1\ \dots \ C_n \big]^t$ in $\C_n(T(\N))$,
 there is an invertible matrix $W \in \fM_n(\T(\N))$ with first row $R$
 and the first column of $W^{-1}$ is $C$.
\item Every atom of $\N$ is finite dimensional, and every non-empty interval
 of $\N$ contains an atom.
\end{enumerate}
\end{thm}

\begin{proof}
Again the discussion prior to the previous theorem shows that $\N$ must
have only finite rank atoms; and there are no non-atomic intervals, namely
every interval contains an atom.  Indeed, it is clear that any interval which
is not finite dimensional must contain infinitely many atoms.
So we need to establish the result for this class of nests.

By Theorem~\ref{T:extend}, the issue is whether the idempotent $P=CR$
is similar in $\fM_n(\T(\N))$ to $E_{11}$.  
We identify $\fM_n(\T(\N))$ with $\T(\N^{(n)})$ acting on $\H^{(n)}$.
By the result of Larson and Pitts~\cite{LP}, it suffices to check that
\[ \rank(P_{E^{(n)}} P P_{E^{(n)}}) = \rank(P_{E^{(n)}} E_{11} P_{E^{(n)}}) \]
and
\[ \rank(P_{E^{(n)}} (I_n-P) P_{E^{(n)}}) = \rank(P_{E^{(n)}} (I_n-E_{11}) P_{E^{(n)}}) \] 
for every interval ${E^{(n)}}$ of $\N^{(n)}$.
When ${E^{(n)}}$ is finite dimensional, it is clear that $P_E R P_{E^{(n)}}$ has
rank $\rank P_E$, and thus this is the rank of the compression of $P$.
This coincides with the rank of the compression of $E_{11}$.
Thus the complements have the complementary ranks, so also agree.
When ${E^{(n)}}$ is infinite dimensional, it contains infinitely many atoms.
The compression of $P$, $I_n-P$, $E_{11}$ and $I_n-E_{11}$ to this interval
will all have non-zero rank in every atom---and so will have infinite rank.
Hence the Larson--Pitts result shows that $P$ is similar to $E_{11}$.
\end{proof}

\section{Split extensions}\label{S:split}

Rieffel \cite{R} establishes several results about ideals and quotients.
One result is that if $\fJ$ is an ideal of $\fA$ with an approximate identity,
then $\rtsr(\fJ) \le \rtsr(\fA)$.  For nonself-adjoint operator algebras such
as nest algebras, the class of ideals with an approximate identity is
extremely limited, and in particular, the only proper ideal in a nest algebra 
with this property is the ideal of all compact operators in the algebra
(see the proof of \cite[Theorem~1.3]{AK}).
Rieffel also shows that $\rtsr(\fA) \le \max\{ \rtsr(\fJ), \rtsr(\fA/\fJ) + 1 \}$.

In this section, we use different hypotheses to obtain similar or better estimates.
We will assume that the quotient map splits.  This will be the case for us
in the application in which the quotient map is the expectation of a nest algebra
into the (block) diagonal corresponding to a subnest.
Among other consequences, we obtain a simpler proof
of \cite[Theorem 2.11]{DLMR}  (Theorem~ \ref{DLMR2.11} above).

{\samepage
\begin{thm} \label{T:split}
Let $\fA$ be a unital Banach algebra, and let $\fJ$ an ideal of $\fA$
such that the quotient map $\pi$ of $\fA$ onto $\fA/\fJ$ is split. 
That is, $0 \to \fJ \to \fA \to \fA/\fJ \to 0$ is split exact.  Then
\begin{enumerate}
\item $\rtsr(\fJ ) \le \rtsr(\fA)$. 

\item If $\tsr(\fA/\fJ) = 1$, then $\rtsr(\fJ ) = \rtsr(\fA)$. 

\item If $\fA/\fJ$ is completely finite, then 
\[ \rtsr (\fA)= \max\{ \rtsr (\fA/\fJ) , \rtsr(\fJ) \} .\]
\end{enumerate}
\end{thm}
}

\begin{proof}
Let the section of $\fA/\fJ$ into $\fA$ be denoted by $\sigma$.
Then $\beta = \sigma\pi$ is an idempotent endomorphism of $\fA$ 
onto a subalgebra $\fB$ isomorphic to the quotient with kernel $\fJ$.
Write $I$ for the unit of $\fA$, and note that it belongs to $\fB$.
Let $\tilde\fJ = \bC I + \fJ$ be the unitization of $\fJ$.

For (i), we first show that whenever $T_i \in \tilde\fJ$ and $X_i \in \fA$ 
satisfy $\sum_{i=1}^n T_iX_i = I$, we can solve this equation
with $X_i \in \tilde\fJ$.
To this end, write $T_i = a_i I + J_i$, where $J_i\in\fJ$,
and $X_i = X^0_i + Y_i$, where $X^0_i=\beta(X_i)$ and $Y_i \in \fJ$.
Then our identity yields
\begin{align*}
 I &= \sum_{i=1}^n T_iX_i  
     = \sum_{i=1}^n a_iX^0_i + \sum_{i=1}^n J_iX^0_i + \sum_{i=1}^n T_i Y_i .
\end{align*}
As the first sum is in $\fB$ and second two sums lie in $\fJ$, it follows that
$\sum_{i=1}^n a_iX^0_i = I$ and $\sum_{i=1}^n J_iX^0_i = - \sum_{i=1}^n T_i Y_i .$

Observe that at least one $a_i \ne 0$.
So one may choose $b_i \in \bC$ so that $\sum_{i=1}^n a_ib_i = 1$.
Set $C_i = b_i I + Y_i + Z_i$ for some unknown $Z_i\in\fJ$.  
Then
\begin{align*}
 \sum_{i=1}^n T_i C_i &= \sum_{i=1}^n (a_i+J_i)(b_i+Y_i+Z_i) \\
 &= \sum_{i=1}^n a_ib_i I  + \sum_{i=1}^n J_i b_i
  + \sum_{i=1}^n T_i Y_i + \sum_{i=1}^n T_iZ_i \\
 &= I  + \sum_{i=1}^n J_i b_i
  - \sum_{i=1}^n J_iX^0_i + \sum_{i=1}^n T_iZ_i 
\end{align*}
Therefore it suffices to solve
\begin{align*}
 \sum_{i=1}^n T_i Z_i &=  \sum_{j=1}^n J_j (X^0_j - b_j) 
 = \sum_{i=1}^n T_iX_i  \sum_{j=1}^n J_j (X^0_j - b_j) \\
 &= \sum_{i=1}^n T_i  \sum_{j=1}^n X_i J_j (X^0_j - b_j) .
\end{align*}
Evidently $Z_i = \sum_{j=1}^n X_i J_j (X^0_j - b_j) \in \fJ$ does the job.

Now suppose that $\rtsr(\fA)=n$, and $T_i = a_iI + J_i$ are given with $J_i\in\fJ$,
and $0 <\ep < 1/2$.  
By changing each $a_i$ if necessary by at most $\ep/2$, we may  and do suppose 
that $|a_i| \ge \ep/2$.
Let $M = \max\{1, \|T_i\| : 1 \le i \le n \}$.
Use  $\rtsr(\fA)=n$ to find $A_i \in \fA$ with $\|A_i\|< \ep^2/6M$ and $X_i \in \fA_i$
so that $\sum_{i=1}^n (T_i+A_i)X_i = I$.
Observe that
\[ T_i+A_i = a_i( I + a_i^{-1}J_i (I + a_i^{-1}A_i)^{-1})) (I + a_i^{-1}A_i) .\]
Set $T'_i = a_i( I + a_i^{-1}J_i (I + a_i^{-1}A_i)^{-1})) = (T_i+A_i) (I + a_i^{-1}A_i)^{-1}$
and note that 
\begin{align*}
 \|T'_i-T_i\| &\le \|T_i\|\, \| (I + a_i^{-1}A_i)^{-1} - I \| + \|A_i\|\, \|(I + a_i^{-1}A_i)^{-1}\| \\
 &\le \frac{ M \| a_i^{-1} A_i \| + \|A_i\|}{1 - \|a_i^{-1} A_i \|}
 \le \frac{ (M \frac2\ep + 1) \frac{\ep^2}{6M}}{ 1 -  \frac2\ep \frac{\ep^2}{6M}} \\
 &\le \frac{\ep}{3} \big(1+\frac{\ep}{6M} \big) \big(1+\frac{\ep}{2M} \big)
  < \frac \ep 2 .
\end{align*}
Then $\sum_{i=1}^n T'_i (I + a_i^{-1}A_i)X_i = I$.
Applying the result of the previous paragraph shows that
there are $Y_i \in \tilde\fJ$ such that $\sum_{i=1}^n T'_i Y_i = I$.
This shows that $\rtsr(\fJ) \le \rtsr(\fA)$.

For (ii), suppose that $\tsr(\fA/\fJ)=1$, and let $n = \rtsr(\fJ)$.
Let $R = \big[ A_1\ \dots \ A_n \big]$ be in $\R_n(\fA)$ and let $\ep > 0$ be given.
Since the invertibles are dense in $\fB$, we may choose
invertible $B_i \in \fB$ so that $\| \beta(A_i) - B_i\| < \ep/2$ for $1 \le i\le n$.
Let $b = \max\{ \|B_i\| : 1 \le i \le n \}$.
Set $A'_i = B_i + (\id-\beta)(A_i)$; so that $\|A_i - A'_i\| < \ep/2$.
Now $A'_i B_i^{-1}$ belong to $\fJ$.

As $\rtsr(\fJ) = n$, there are elements $J_i \in \tilde\fJ$ with $\|J_i\| < \frac{\ep}{2b}$
and $X_i \in \tilde\fJ$ such that
\[
 I = \sum_{i=1}^n (A'_iB_i^{-1} + J_i ) X_i 
   = \sum_{i=1}^n ( A'_i + J_iB_i) (B_i^{-1} X_i)
\]
Now $\|A_i - ( A'_i + J_iB_i)\| < \frac\ep 2 + \|J_i\|\,\|B_i\| < \ep$.
Thus $\rtsr(\fA) \le \rtsr(\fJ)$.
The reverse inequality follows from (i).

For (iii), we assume that $\fB$ is completely finite;
and let $n =  \max\{ \rtsr (\fA/\fJ) , \rtsr(\fJ) \}$.  We are given
$R = \big[ A_1\ \dots \ A_n \big] \in \R_n(\fA)$ and $\ep > 0$.
Because $\beta(R) = \big[ \beta(A_1)\ \dots \ \beta(A_n) \big]$ is in $\R_n(\fB)$ 
and $\rtsr(\fB) \le n$, there is an $\ep/2$ perturbation
$S = \big[ B_1\ \dots \ B_n \big]$ with right inverse $C \in \C_n(\fB)$.
As $\fB$ is completely finite, there is an invertible matrix $W \in \fM_n(\fB)$
with first column equal to $C$.
Then $SW =  \big[ I \ 0 \ \dots \ 0 \big]$.

Let $A'_i = B_i + (\id-\beta)(A_i)$ and $R' =  \big[ A'_1\ \dots \ A'_n \big]$.
It follows that  $\|R-R'\| < \ep/2$, and that $R'W =  \big[ I+J_1\ J_2\ \dots \ J_n \big]$
lies in $\R_n(\tilde\fJ)$.
Since $\rtsr(\fJ) \le n$, there is a perturbation $K \in \R_n(\tilde\fJ)$
with $\|K\|  \le \frac \ep{2\|W^{-1}\|}$
and $X \in \C_n(\tilde\fJ)$ so that 
\[
 I = (R'W+K)X = (R'+KW^{-1})(WX) .
\]
As $\|R-(R'+KW^{-1})\| < \ep/2 + \|K\|\,\|W^{-1}\| < \ep$, we obtain
$\rtsr(\fA) \le n$.

The reverse inequality follows for $\fJ$ by (i), 
and for $\fA/\fJ$ by \cite[Theorem~4.3]{R}.
\end{proof}

If $\fA_k$ are unital operator algebras acting on Hilbert spaces $\H_k$,
let $\bigoplus_{k=1}^\infty \fA_n$ denote the $\linf$ direct sum
(or direct product)  of all operators of the form 
$A = \bigoplus_{k=1}^\infty A_k$ acting on $\bigoplus_{k=1}^\infty \H_k$
such that $A_k \in \fA_k$ and $\sup_{k\ge1} \|A_k\| < \infty$.

An immediate consequence of this for nests is the following.

\begin{lem}\label{L:split nest}
Let $\N$ be a nest of order type $\omega$ with finite rank atoms. 
Let $\M$ be an infinite subnest with atoms $E_i$.
Let $\Delta_\M$ denote the expectation onto the diagonal $\fD$ of $\T(\M)$;
and let $\fB = \Delta_\M(\T(\N))$.  Then
\[ \rtsr(\T(\N)) = \max\{ \rtsr(\T(\M)), \rtsr(\fB) \} .\]
\end{lem}
 
\begin{proof}
Let $\fJ = \ker \Delta_\M$.
Since $\fD$ is a finite type I von Neumann algebra, $\tsr(\fD)=1$.
Now $\fJ$ is also an ideal in $\T(\N)$ and 
$\T(\N)/\fJ \simeq \fB$.
Both of these quotient maps split by the natural identification
of $\fD$ and $\fB$ as subalgebras of $\T(\M)$ and $\T(\N)$, respectively.

Hence $\rtsr(\fJ) = \rtsr(\T(\M))$ by Theorem~\ref{T:split}.
Using Theorem~\ref{T:split} again, we obtain 
\[ 
 \rtsr(\T(\N)) =  \max\{ \rtsr(\fJ), \rtsr(\fB) \} 
 =  \max\{ \rtsr(\T(\M)), \rtsr(\fB) \} .
\qedhere
\]
\end{proof}

\begin{thm} \label{T:chain of nests}
Let $\N$ be a nest of order type $\omega$ with finite rank atoms.
Suppose that a sequence of infinite subnests 
$\N = \N_0 \supset \N_1 \supset \dots \supset \N_k$ is given.
Then 
\[
 \rtsr(\T(\N)) = \max\big\{ 
 \rtsr(\T(\N_k)), \rtsr(\Delta_{\N_j}(\T(\N_{j-1}))), 1 \le j \le k \big\} .
\]
\end{thm}

\begin{proof}
Let $\fB_{j-1} = \Delta_{\N_j}(\T(\N_{j-1}))$ for $1 \le j \le k$.
This is just a repeated application of Lemma~\ref{L:split nest}.
We obtain 
\[ \rtsr(\T(\N_{k-1})) = \max\{ \rtsr(\T(\N_k)), \rtsr(\fB_{k-1}) \} .\]
Then 
\begin{align*}
 \rtsr(\T(\N_{k-2})) &= \max\{ \rtsr(\T(\N_{k-1})), \rtsr(\fB_{k-2}) \} \\
  &= \max\{ \rtsr(\T(\N_k)), \rtsr(\fB_{k-1}), \rtsr(\fB_{k-2}) \}
\end{align*}
After $k$ steps, we arrive at the desired conclusion.
\end{proof}

\begin{defn}
Consider a nest $\N=\{N_i, \H : i \ge0\}$ of order type $\omega$ 
and an infinite subnest $\M = \{N_{k_i},\H : i \ge 0\}$.
Say that $\M$ has \textit{finite index} in $\N$ 
if $\sup_{i\ge0} k_{i+1}-k_i = l < \infty$.
\end{defn}

\begin{cor} \label{C:finite index}
Suppose that $\N=\{N_i, \H : i \ge0\}$ is a nest of order type $\omega$
with finite rank atoms.  Let $\M$ be a subnest of finite index in $\N$.
Then $\rtsr(\T(\M))=\rtsr(\T(\N))$.
\end{cor}

\begin{proof}
Let $E_i = N_{k_i}\ominus N_{k_{i-1}}$ for $i \ge 1$ be the atoms of $\M$.
Then we observe that $\fB = \Delta_\M(\T(\N)) = \bigoplus_{i\ge1} \T(E_i \cap \N)$.

Each $\T(E_i \cap \N)$ is a nest algebra on a finite dimensional space
with nest of length at most $l$.
An element $T$ in this algebra has a block matrix $\big[ T_{ij} \big]$
where $T_{ij}=0$ when $i>j$.
Thus $T$ is invertible if and only if each $T_{ii}$ invertible.
Moreover, since each $T_{ii}$ belongs to a full matrix algebra, polar decomposition
allows one to find an $\ep$ perturbation with inverse bounded by $\ep^{-1}$.
The norm of the inverse is bounded by $C_l \|T\|^{-1} \ep^{-l}$ where 
$C_l$ is a constant.
See Remark~\ref{R:finite nests} in the next section for more detail.
Consequently we conclude that the invertibles are dense in $\fB$.
That is, $\tsr(\fB) = 1$.
Alternatively, one can observe that $\fB$ may be considered as having an
$l\times l$ upper triangular form with finite von Neumann algebras in the
diagonal entries.  The invertibles are dense in these diagonal entries, and
when the diagonal is invertible, so it the whole operator. 

The conclusion now follows from Corollary~\ref{L:split nest}.
\end{proof}

This leads to a simplification of the proof of Theorem~\ref{DLMR2.11}.
The hypothesis is that
\[  \max_{1 \le i \le  (k+1)J} n_i \ge r  \max_{1 \le i \le  kJ} n_i \qforal k \ge 1 \]
for some integer $J\ge1$ and real number $r>1$.
Choose an integer $p \ge 1$ so that $r^p \ge 5 p J$.  
Then take $k_i = pJi$ for $i\ge0$; and let $\M=\{N_{k_i},\H : i \ge1\}$.
An easy argument from the proof of \cite[Theorem 2.11]{DLMR} shows that  
$\dim(N_{k_i}\ominus N_{k_{i-1}}) \ge 4 \dim N_{k_{i-1}}$ for each $i \ge 2$.	
By Corollary~\ref{C:finite index}, $\rtsr(\T(\N))=\rtsr(\T(\M))$.
Now $\rtsr(\T(\M)) = 2$ follows from \cite[Theorem~2.1]{DLMR}, which has an easier proof.

Another immediate consequence of Lemma~\ref{L:split nest} is the following useful fact.

\begin{cor} \label{L:subnest}
Let $\N$ be a nest of order type $\omega$ with finite rank atoms, 
and let $\M$ be an infinite subnest of $\N$.
Then \[ \rtsr(\T(\M)) \le \rtsr(\T(\N)) .\]
\end{cor}

\section{Quantitative measurement of tsr}\label{S:quant}

It will be useful to have a quantitative measurement of the stable rank
in order to establish an obstruction.

\begin{defn}
Given a Banach algebra $\fA$, an element $A \in \R_n(\fA)$ and $t>0$, let
\[
 \rho_n(A,t) \!=\! 
 \inf \{  \|C\| :\! C \in \C_n(\fA), \exists B \!\in\! \R_n(\fA), \|A-B\| < t\, \text{and}\, BC = I \} .
\]
By convention, the $\inf$ over the empty set is $+\infty$. 
Then let the \textit{right stable range function} on $\fA$ be
\[
 \rho_n(\fA,t) = 
 \sup \{ \rho_n(A,t) : A \in \R_n(\fA),\ \|A\| \le 1 \} 
 \qfor 0<t \le 1 .
\]
\end{defn}

Similarly,we can define the \textit{left stable range function} $ \lambda_n(\fA,t)$ on $\fA$.
We will only deal with the right version in this paper; but clearly there is always a 
corresponding left version of each theorem.
There is no point in considering $t>1$; and indeed, we are only interested in 
the behaviour as $t\to 0$.

One simple property of this function is the relation to quotients.

\begin{prop}\label{P:quotients}
Let $\fJ$ be an ideal of $\fA$.  
Then $\rho_n(\fA/\fJ,t) \le \rho_n(\fA,t)$ for all $n \ge 1$ and $0<t \le 1$.
\end{prop}

\begin{proof}
Let $\pi$ be the quotient map of $\fA$ onto $\fA/\fJ$.
Given $\dot A \in \R_n(\fA/\fJ)$ with $\|\dot A \| < 1$,
select an element $A \in \R_n(\fA)$ with $\pi(A) = \dot A$ and $\|A\|<1$.
Then there is a perturbation $B \in \R_n(\fA)$ with $\|B-A\|<t$ and
an element $C \in \C_n(\fA)$ so that $BC=I$ and $\|C\| \le \rho_n(\fA,t)$.
Therefore $\dot B = \pi(B)$ and $\dot C = \pi (C)$ show that
$\rho_n(\dot A,t) \le  \rho_n(\fA,t)$.
\end{proof}

For $\rtsr(\fA) \le n$, one must have $\rho_n(A,t) < \infty$ 
for all $A \in \R_n(\fA)$ and all $0<t \le 1$.
In all cases that we understand, this is established by showing 
that $\rho_n(\fA,t) < \infty$ for all $0<t \le 1$.

\begin{ques}
If $\fA$ is an Banach algebra with $\rtsr(\fA) =n < \infty$, 
is $\rho_n(\fA,t) < \infty$ for all $0<t \le 1$?
\end{ques}

We are able to establish this in the case of interest to us.

\begin{prop} \label{P:direct sum}
Let $\fA_k$ be unital operator algebras for $k\ge1$, 
and let $\fA = \bigoplus_{k=1}^\infty \fA_n$.  Then
\[ \rho_n(\fA,t) = \sup_{k\ge1} \rho_n(\fA_k,t) .\]
Moreover, if $\rho_n(\fA_k,t)< \infty$ for all $k \ge 1$ and all $0<t \le 1$, then
$\rtsr(\fA) \le n$ if and only if $ \rho_n(\fA,t)< \infty$ for all $0<t \le 1$.
\end{prop}

\begin{proof}
The first claim follows from the elementary observation that
\[ \rho_n(\textstyle\bigoplus_{k=1}^\infty A_k, t) = \sup_{k\ge1} \rho_n(A_k,t) .\]
In particular, if each $\rho_n(\fA_k,t)< \infty$ but $\rho_n(\fA,t) = \infty$, then
there must be $A_k \in \R_n(\fA_k)$
with $\|A_k\|\le 1$ and $\sup_{k\ge1} \rho_n(A_k,t) = \infty$.
Thus $\rho_n(\textstyle\bigoplus_{k=1}^\infty A_k, t) = \infty$ and so $\rtsr(\fA) > n$.
The other direction is easy.
\end{proof}

Theorem~\ref{T:split} can be quantified by examining the details of the proof.
The exact relationships are not so important.  What we require is some
control based on the data.
Part (iii) uses Theorem~\ref{T:extend}.
In order to obtain a quantitative version, control on the norm of the
similarity is required. So we make the following definition.

\begin{defn}
If $\fA$ is completely finite, say that it is \textit{uniformly completely finite}
if there are real functions $\sigma_n$ so that whenever $R$ in $\R_n(\fA)$ and
$C$ in $\C_n(\fA)$ satisfy $RC = 1$, then the invertible operator $W$ with
first row $R$ and inverse $W^{-1}$ with first column $C$ satisfy
\[ \|W\| \le \|R\| \sigma_n(\|R\|\,\|C\|) \qand  \|W^{-1}\| \le \|C\| \sigma_n(\|R\|\,\|C\|) .\]
\end{defn}

We will show below that completely finite nest algebras are uniformly completely finite
using $\sigma_n(x) = 2 x$.
In our applications to nest algebras,  the expectation $\beta$ is completely contractive;
so $\| \beta \otimes \id_n\| =1$ for all $n \ge 1$.

\pagebreak[3]
{\samepage
\begin{thm}  \label{T:quantitative split}
Let that $\fA$ be a unital Banach algebra, and let $\fJ$ an ideal of $\fA$
such that the quotient map $\pi$ of $\fA$ onto $\fA/\fJ$ is split. 
Then
\begin{enumerate}
\item $\rho_n(\fJ,t) \le 7  \| \beta \otimes \id_n\|^2 \rho_n(\fA, \frac{t^2}6)^2$. 

\item If $\tsr(\fA/\fJ) = 1$, then $\rho_n(\fA,t ) \le 
\rho_1(\fB,\frac t 2) \, \rho_n \big( \fJ,\frac t{2\|\beta\| \rho_1(\fB,t/2)} \big)$. 

\item If $\fA/\fJ$ is uniformly completely finite with respect to functions $\sigma_n$,
then 
\[
 \rho_n (\fA,t) \le \sigma_n\big( \rho_n(\fB,\tfrac t 2)\big)  
 \rho_n \big( \fJ, \tfrac t{2\rho_n(\fB,t/2)\,\sigma_n\rho_n(\fB,t/2)} \big) .
\]
\end{enumerate}
\end{thm}
}

\begin{proof}
We only provide a sketch.  Let $B_n = \| \beta \otimes \id_n\|$.  

First suppose that $T \in \R_n(\tilde\fJ)$ has $\|T\| \le 1$ and $X \in \C_n(\fA)$
satisfies $TX = I$.  Then we found $Y \in \C_n(\tilde\fJ)$ such that $TY = I$.\vspace{.3ex}
If we also assume that $\max\{ |a_i| : 1 \le i \le n \} \ge t$, then we can
choose $\sum a_ib_i = 1$ with $\sum |b_i| \le t^{-1}$.
Then $Y = \big[ Y_1\ \dots\ Y_n\big]^t$ is given by
\[ Y_i = b_i + (\id - \beta)(X_i) + X_i \sum_{j=1}^n (\id-\beta)(T_j) (\beta(X_j)-b_j) .\]
Thus we obtain $\|Y\| \le (B_n^2+B_n)\|X\|^2 + (2B_n+3)t^{-1}\|X\|$.

Then following through the second calculation, we fixed $t>0$, 
made a $t/2$ perturbation to ensure that $|a_i| \ge t/2$ for each $i$,
and then found a $t^2/6$ perturbation to obtain a solution 
$\sum (T_i+A_i)X_i = I$ in $\fA$.  So we may take $\|X\| \le \rho_n(\fA, t^2/6)$.
Also $\rho_n(\fA,t) \ge t^{-1}$.
Therefore, by the previous paragraph,  the solution $T'Y=I$ in $\tilde\fJ$
satisfies 
\[ \|Y\| \le (B_n^2+B_n)\|X\|^2 + (2B_n+3)t^{-1}\|X\| \le 7B_n^2 \rho_n(\fA, t^2/6)^2 .\] 

For (ii), we note that $\|B_i^{-1}\|$ can be chosen to be at most $\rho_1(\fB,t/2)$.
Then note that if $\|A\| = M$, then there is a perturbation $A'$ of $A$ of norm at most $t>0$
so that $A'X=I$ has a solution with $\|A' - A\| \le M \rho_n(\fA, t/M)$.
In our case, we are taking a $t/(2\|\beta\|)$ perturbation, which leads to the
given estimate.

Finally consider (iii).  Starting with $R \in \R_n(\fA)$ with $\|R\|\le 1$ and $t>0$,
we first find a $t/2$ perturbation of $\beta\otimes\id_n(R)$ with right inverse
$C \in \C_n(\fB)$ with $\|C\|  \le \rho_n(\fB,t/2)$.
Then one finds an invertible operator $W$ with
first row $R$ and inverse $W^{-1}$ with first column $C$ such that
$\|W\| \le \sigma_n(\|C\|)$ and  $\|W^{-1}\| \le \|C\| \sigma_n(\|C\|)$.
Working through the rest of the details, one obtains $R'W \in \R_n(\tilde\fJ)$
and a perturbation $K$ of norm at most $(2\|C\| \sigma_n(\|C\|))^{-1}t$ 
so that $R'W+K$ has a right inverse $X$ in $\C_n(\tilde\fJ)$ 
with $\|X\| \le \rho_n \big( \fJ, \tfrac t{2\rho_n(\fB,t/2)\,\sigma_n\rho_n(\fB,t/2)} \big)$. 
Finally the right inverse of the final $t$-perturbation of $R$ is $WX$.
So
\[ 
 \rho_n (\fA,t) \le \sigma_n\big( \rho_n(\fB,\tfrac t 2)\big)  
 \rho_n \big( \fJ, \tfrac t{2\rho_n(\fB,t/2)\,\sigma_n\rho_n(\fB,t/2)} \big) .
 \qedhere
\]
\end{proof}

Here is the necessary estimate for nest algebras.

\begin{prop} \label{P:uniform comp. finite}
A completely finite nest algebra is uniformly completely finite with
respect to the functions $\sigma_n(x) = \sigma(x) =2 x$.
\end{prop}

\begin{proof}
The proof of Theorem~\ref{T:extend_nest} contains an explicit description of the 
similarity between $P=CR$ and $E_{11}$, namely $UV$ where $U$
is unitary and $Vx = Px \oplus (I-P) x$.
So 
\[
 \|V\| \le \big( \|P\|^2 + \|I-P\|^2 \big)^{1/2} 
 = \sqrt 2 \|P\| \le \sqrt 2 \|R\|\,\|C\| .
\]
Since $V^{-1}(x \oplus y) = x+y$, it follows that $\|V^{-1}\| \le \sqrt2$
and therefore 
\[ \|V\|\,\|V^{-1}\| \le 2 \|R\|\,\|C\| .\]

Chasing through the calculation in the proof of Theorem~\ref{T:extend}
to show that (iv) implies (iii), one finds that
\[
 \|W\| \le \|R\|\,\|V^{-1}\|\, \|V\| \le 2 \|R\|^2\,\|C\|  
\]
and 
\[ 
 \|W^{-1}\| \le \|V^{-1}\|\,\|V\|\,\|C\| \le 2 \|R\|\,\|C\|^2 . \qedhere
\]
\end{proof}

Here is a simple but useful estimate.

\begin{lem}\label{L:2x2}
Suppose that $\fX$ is an operator $\fA_1$--$\fA_2$ bimodule, and form the
matrix algebra $\fA = \begin{bmatrix}\fA_1 & \fX\\0 & \fA_2\end{bmatrix}$.
Then $\rtsr(\fA) = \max\{ \rtsr(\fA_1), \rtsr(\fA_2) \}$ and 
\begin{align*}
 \max\{\rho_n(\fA_1,t), \rho_n(\fA_2,t)\} \le \rho_n(\fA,t)  &\le 
 \max\{\rho_n(\fA_1,t), \rho_n(\fA_2,t)\} + \rho_n(\fA_1,t) \rho_n(\fA_2,t) .
\end{align*}
\end{lem}

\begin{proof}
Take any $A \in \R_n(\fA)$ with $\|A\|<1$. \vspace{.4ex}
We may write this as $A = \begin{bmatrix}A_{11}&A_{12}\\0&A_{22}\end{bmatrix}$
where $A_{11} \in \R_n(\fA_1)$, $A_{22} \in \R_n(\fA_2)$
and $A_{12} \in \R_n(\fX)$.
Choose perturbations $B_{11} \in \R_n(\fA_1)$ and $B_{22} \in \R_n(\fA_2)$ 
such that $\|A_{11}-B_{11}\| < t$ and $\|A_{22} - B_{22} \| < t$, 
and elements $C_{11} \in \C_n(\fA_1)$ and $C_{22} \in \C_n(\fA_2)$ 
with $\|C_{11}\| \le \rho_n(\fA_1,t)$ and $\|C_{22}\| \le \rho_n(\fA_2,t)$ 
so that $B_{11}C_{11} = I_{\fA_1}$ and $B_{22}C_{22} = I_{\fA_2}$.
Then define
\[
 B = \begin{bmatrix}B_{11}&A_{12}\\0&B_{22}\end{bmatrix} \qand 
 C = \begin{bmatrix}C_{11}& -C_{11} A_{12} C_{22} \\0&C_{22}\end{bmatrix}.
\]
One readily verifies that $\|A-B\| < t$ and $BC = I_\fA$.
Moreover we obtain the simple estimate
\begin{align*}
 \|C\| &\le \max\{\|C_{11}\|, \|C_{22}\| \} + \|C_{11}\|\,\|C_{22}\| \\
  &\le \max\{\rho_n(\fA_1,t), \rho_n(\fA_2,t)\} + \rho_n(\fA_1,t) \rho_n(\fA_2,t) .
\end{align*}
This establishes the upper bound for $\rho_n(\fA,t)$.  The lower bound is evident.
\end{proof}

\begin{rem} \label{R:finite nests}
Since $\fM_k$ is a finite dimensional C*-algebra,
a simple use of the polar decomposition shows that $\rho_n(\fM_k,t) = 1/t$.

Now suppose that $\N = \{N_0 < N_1 < \dots < N_p = \H\}$ is a nest on a 
\textit{finite dimensional} Hilbert space $\H$.
Then $\T(\N)$ has a $p\times p$ block upper triangular form.
The diagonal entries of $\T(\N)$ are full matrix algebras
acting on finite dimensional spaces.
Hence a repeated application of Lemma~\ref{L:2x2} shows that 
$\rho_n(\T(\N),t) = O(t^{-p})$ independent of $\dim\H$.
\end{rem}

\begin{thm}\label{T:sup_rtsr}
Let $\N = \{ N_k, \H : k \ge 0 \}$ be a nest of order type 
$\omega$ with finite dimensional atoms.
Let $\N_k = \{N_0, N_1,\dots,N_k\}$ be the restriction of $\N$ to $N_k$.
Then 
\[ \rho_n(\T(\N),t) = \sup_{k\ge1} \rho_n(\T(\N_k),t) ;\] 
and $\rtsr(\T(\N)) \le n$ if and only if  
$\rho_n(\T(\N),t) < \infty$ for all $0<t \le 1$.
\end{thm}

\begin{proof}
Clearly, $R := \sup_{k\ge1} \rho_n(\T(\N_k),t) \le \rho_n(\T(\N),t)$.
Also, since $\T(\N_k)$ are finite dimensional, it is easy to see that the
invertible elements are dense---and so $\tsr(\T(\N_k)) = 1$ for all $k\ge1$.
By the previous Remark, we know that $\rho_n(\T(\N_k),t) < \infty$ for all
$n,k \ge 1$ and all $0<t \le 1$.

On the other hand, suppose that $A \in \R_n(\T(\N))$ with $\|A\|< 1$.
Then for each $k\ge 1$, consider $A_k := A P_{N_k} \in \R_n(\T(\N_k))$.
We can find $B_k \in \R_n(\T(\N_k))$ and $C_k \in\C_n(\T(\N_k))$ 
such that $\|A_k-B_k\| < t$,
$\|C_k\| \le R$ and $B_k C_k = P_{N_k}$.
Consider $B_k$ as elements of $\R_n(\T(\N))$ which vanish on $N_k^\perp$
and likewise for $C_k$.
Choose a subnet of this sequence so that both 
\[
 \textstyle
 \wotlim_\mu B_{k(\mu)} = B 
 \qand  
 \wotlim_\mu C_{k(\mu)} = C
\]
exist.
Since $\T(\N)$ is \wot-closed, these limits lie in $\R_n(\T(\N))$ and
$\C_n(\T(\N))$ respectively; and $\|A-B\| \le t$ and $\|C\| \le R$.
Multiplication is jointly \wot-continuous on $\T(\N)$ for the simple reason that
the computation of the $i,j$ entry of the product of two upper triangular
operators is a finite sum of products of operators on finite dimensional space.
Therefore $BC = I$. Hence $\rho_n(\T(\N),t) \le R$.

Next suppose that the supremum is infinite.
We claim that 
\[ \sup_{k>k_0} \rho_n(\T((N_k \ominus N_{k_0})\cap \N),t) = \infty \]
for any $k_0$.
Indeed, we can write $\T(N_k \cap \N)$ in the matrix form
\[
 \begin{bmatrix}
 \T(N_{k_0} \cap \N) & \B(N_k \ominus N_{k_0}, N_{k_0}) \\
 0 & \T((N_k \ominus N_{k_0})\cap \N) 
 \end{bmatrix} .
\]
So by Lemma~\ref{L:2x2}, 
$ \rho_n(\T(N_k\cap \N),t) \le \max\{X_0,Y_k\}+X_0Y_k $
where $X_0 = \rho_n(\T(N_{k_0}\cap \N),t)$ and 
$Y_k = \rho_n(\T(N_k \ominus N_{k_0})\cap \N),t)$.
Since the left hand side is unbounded, so must the $Y_k$'s be unbounded.
We conclude that we can recursively select a sequence $k_j$ so that 
$E_j = N_{k_j} \ominus N_{k_{j-1}}$ satisfy $\rho_n(\T(E_j\cap \N),t) \ge j$.

There is a compression homomorphism
\[ \pi(T) = \textstyle\bigoplus_{j=1}^\infty P_{E_j}T|_{E_j} \]
onto $\fB = \bigoplus_{j=1}^\infty \T(E_j\cap \N)$.
By Proposition~\ref{P:direct sum}, this algebra has $\rtsr(\fB) > n$.
But homomorphic images have smaller stable rank.  So $\rtsr(\T(\N)) > n$ also.
\end{proof}

\section{Upper and Lower bounds for nest algebras} \label{S:bounds}

We now apply the results of the previous section to nest algebras.
The goal is to obtain quantitative estimates that will allow us to
show that the right topological stable rank is either at most 2
or it is infinity.

\cite[Theorem~2.1]{DLMR} shows that if $\N$ is an atomic nest order 
isomorphic to $\omega$ with finite dimensional atoms of rank $n_k$ \vspace{.2ex}
satisfying $n_k \ge 4 \sum_{i<k} n_i$, then $\rtsr(\T(\N)) =2$.
The proof is an explicit construction, and if one checks the proof,
one easily obtains the following quantitative estimate:

\begin{lem} \label{L:rapid growth}
Let $\N$ be an atomic nest order  isomorphic to $\omega$ with 
finite dimensional atoms of  rank $n_k$ 
satisfying $n_k \ge 4 \sum_{i<k} n_i$. 
Then
\[ \rho_2(\T(\N)),t) < 9 t^{-2} .\]
\end{lem}

If we now apply Theorem~\ref{T:quantitative split}, we can obtain a
quantitative version of Theorem~\ref{DLMR2.11}.

\pagebreak[3]
{\samepage 
\begin{thm}\label{T:quantitative 2.11}
Suppose that $\N$ is a nest of order type $\omega$ with finite rank
atoms of dimensions $n_i$ such that for some integer $J\ge1$ and 
real number $r>1$,
\[ \max_{i \le (k+1)J} n_i \ge r  \max_{i \le kJ} n_i \qforal k \ge 1. \]
If an integer $p$ satisfies $r^p>5pJ$, 
then 
\[ \rho_2(\T(\N),t) \le C (t^{-5pJ-4}) \]
where $C$ is a constant depending only on $pJ$.
\end{thm}
}

\begin{proof}
Take $k_i = pJi$ for $i\ge0$, and let $\M=\{N_{k_i},\H : i \ge 0\}$.
An easy argument from the proof of \cite[Theorem 2.11]{DLMR} shows that  
$\dim(N_{k_{i+1}}\ominus N_{k_i}) \ge 4 \dim N_{k_i}$ for each $i \ge 1$.	
Hence Lemma~\ref{L:rapid growth} applies to $\M$, so that
$\rho_2(\T(\M),t) < 9t^{-2}$.

Let $\Delta_\M$ be the expectation onto the diagonal of $\T(\M)$,
and let $\fJ = \ker \Delta_\M$.
By Theorem~\ref{T:quantitative split}, 
$\rho_2(\fJ,t) < 7\big(9(t^2/6)^{-2}\big)^2 = C t^{-4}$.
Now $\fJ$ is also an ideal in $\T(\N)$, and 
$\T(\N)/\fJ \simeq \fB = \Delta_\M(\T(\N))$.
Since $\M$ has finite index $pJ$ in $\N$, $\tsr(\fB) = 1$.
Indeed, $\fB$ is a direct sum of $pJ \times pJ$ block upper triangular
matrix algebras.
By Remark~\ref{R:finite nests}, $\rho_1(\fB,t) < C' t^{-pJ}$
where the constant $C'$ depends on $pJ$.

Hence using Theorem~\ref{T:quantitative split} again, we obtain
\[ 
  \rho_2(\T(\N),t) \le C' (t/2)^{-pJ} C \big(t(2C' (t/2)^{-pJ})^{-1}\big)^{-4}
  = C'' t^{-5pJ-4} .
  \qedhere
\]
\end{proof}

Now we turn to the problem of obtaining lower bounds, which will allow
us to establish that the right stable rank is infinity in many cases.
The following lemma captures the key idea.

If $\M = \{ \{0\}=M_0 < M_1 < \dots < M_p\}$ and $\N = \{ \{0\}=N_0 < N_1 < \dots < N_p\}$
are two finite nests, let 
\[
 \T(\M,\N) = \{ T \in \B(M_p,N_p) : TM_j \subset N_j \FOR 1 \le j \le p \}
\]
and
\[
 \T_0(\M,\N) = \{ T \in \B(M_p,N_p) : TM_j \subset N_{j-1} \FOR 1 \le j \le p \} .
\]

\begin{lem}\label{L:lower bound}
Suppose that 
\[ \dim M_j \ominus M_{j-1} \le \dim N_{j-1} \ominus N_{j-2} \qfor 2 \le j \le p .\]
Then there is a partial isometry $A \in \T_0(\M,\N)$ such that
for any $B$ in $\T(\M,\N)$ satisfying $\|A-B\| \le t \le \tfrac14$, we have
\[ \inf \big\{ \|C\| : C \in \T(\M,\N) \AND BC=I \big\} \ge (1-2t)^{p-1} t^{-p} \ge 2^p .\]
\end{lem}

\begin{proof}
Let $A_j$ be an isometry of the atom $E_j := M_j \ominus M_{j-1}$ into 
 \mbox{$F_{j-1} := N_{j-1} \ominus N_{j-2}$} for $2 \le j \le p$.
Take $A = \sum_{j=2}^p A_j$.  Let $B$ be any upper triangular perturbation
of norm at most $t$, and write $B$ as a matrix $\big[ B_{ij} \big]$ with respect
to the decomposition $M_p = \bigoplus_{j=1}^p E_j$ and $N_p=\bigoplus_{j=1}^p F_j$.
Then $B_{ij} = 0$ for $j<i$, $\|A_j - B_{j-1,j}\| \le t$ and $\|B_{ij}\| \le t$ in all other cases.
Write any right inverse $C \in \T(\N,\M)$ of $B$ 
as an upper triangular matrix $\big[ C_{ij} \big]$.

Observe that $B_{pp} C_{pp} = P_{E_p}$; and hence $\|C_{pp}\| \ge t^{-1}$.
We will estimate $\|C_{jp}\|$ by induction.
Note that $B_{p-1,p-1}C_{p-1,p} + B_{p-1,p} C_{p,p} = 0$.
As $A_p$ is isometric on $E_p$, the operator $B_{p-1,p}$ is bounded below
by $1-t$.  Hence
\begin{align*}
 \| C_{p-1,p} \| &\ge \|B_{p-1,p} C_{p,p}\| \|B_{p-1,p-1}\|^{-1}  \\
 &\ge (1-t) t^{-1} \|C_{p,p}\| \ge (1-2t) t^{-1} \|C_{p,p}\|.
\end{align*}

Assume that we have shown that $\|C_{j,p}\| \ge (1-2t)t^{-1}\|C_{j+1,p}\|$ for $k<j<p$.
Then $\|C_{j+s,p}\| \le t^s (1-2t)^{-s} \|C_{jp}\|$ for $1 \le s \le p-j$.
We have $\sum_{j=k}^p B_{kj}C_{jp} = 0$.   So compute
\begin{align*}
 \|C_{kp}\| &\ge 
 \|B_{kk}\|^{-1} \Big( \|B_{k,k+1} C_{k+1,p}\| - \sum_{j=k+2}^p \|B_{kj}\|\,\|C_{jp}\| \Big) \\
 &\ge t^{-1} \Big( (1-t) \|C_{k+1,p}\| -  \sum_{s=1}^{p-k-1} t \|C_{k+1+s,p}\| \Big)  \\
 &\ge t^{-1} \Big( (1-t) \|C_{k+1,p}\| -  \sum_{s=1}^{p-k-1} t t^s (1-2t)^{-s} \|C_{k+1,p}\| \Big)  \\
 &\ge t^{-1} \|C_{k+1,p}\| \Big( 1-t- t \sum_{s\ge1} \big( \tfrac{t}{1-2t} \big)^s  \Big) \\ 
 &\ge t^{-1} (1-2t) \|C_{k+1,p}\| .
\end{align*}
The final inequality uses the estimate $t \le 1/4$.
Finally we obtain
\[ 
 \|C\| \ge  \|C_{1p}\| \ge t^{1-p} (1-2t)^{p-1} \|C_{pp}\| \ge (1-2t)^{p-1} t^{-p} \ge 2^p . 
 \qedhere
\]
\end{proof}

This lemma is applied to show that a sufficiently long strings of consecutive
intervals of a nest of \textit{decreasing} dimension lead to a quantitative obstruction
to finite topological stable rank.

{\samepage
\begin{thm} \label{T:lower bound}
Let $\N$ be a nest, and let $q = (n^p-1)/(n-1)$.  
Suppose that $\N$ contains a chain of elements
$N_0 < N_1 < \dots < N_q$ so that
\[ \dim(N_j \ominus N_{j-1}) \ge \dim(N_{j+1} \ominus N_j) \qfor 1 \le j < q .\]
Then 
\[ \rho_n(\T(\N),t) \ge (1-2t)^{p-1} t^{-p} \ge 2^p \qforal 0<t\le 1/4. \]
\end{thm}
}

\begin{proof}
The restriction of $\N$ to the interval $N_q \ominus N_0$ is a contractive homomorphism.
So by Proposition~\ref{P:quotients}, it suffices to establish the result for the restricted nest.
Let \[ k_s = n^{p-1}+\dots + n^{p-s} \qfor 1 \le s \le p .\]
Let $\N' = \{N_0 < N_{k_1} < N_{k_2} < \dots < N_{k_p} \}$.
Observe that the hypotheses ensure that
\[
 \dim(N_{k_i} \ominus N_{k_{i-1}}) \ge n \dim(N_{k_{i+1}} \ominus N_{k_i})
 \qforal 1 \le i < p.
\]
Also note that $\bC I + \T_0(\N') \subset \T(\N) \subset \T(\N')$.

Consider an element of $\R_n(\T(\N'))$ as an element of $\T(\M,\N')$
where 
\[ \M = \N^{\prime(n)} = \{N_0^{(n)} < N_{j_1}^{(n)} <  \dots < N_{j_p}^{(n)} \} .\]
Here $N^{(n)}$ is the direct sum of $n$ copies of $N$.
So the inequalities of the previous paragraph may be restated as
\[
 \dim(N_{k_i} \ominus N_{k_{i-1}}) \ge  \dim(N_{k_{i+1}}^{(n)} \ominus N_{k_i}^{(n)})
 \qforal 1 \le i < p.
\]
Thus Lemma~\ref{L:lower bound} applies.

We obtain a partial isometry $A \in \T_0(\M,\N) \subset \R_n(\T(\N))$
so that 
\[ \rho_n(A,t) \ge (1-2t)^{p-1} t^{-p} \qforal 0<t\le 1/4 .\] 
where this is computed relative to the algebra $\T(\N')$.  But computing it with respect
to the subalgebra $\T(\N)$ can only increase this value.
Hence 
\[  \rho_n(\T(\N),t) \ge (1-2t)^{p-1} t^{-p} .\qedhere \]
\end{proof}

Now we obtain a very useful way to establish that the right stable rank is infinite.

\begin{cor}\label{C:decreasing}
Let $\N$ be a nest.  Suppose that for every positive integer $q$, 
there is a chain $N_0<N_1<\dots < N_q$ in $\N$ so that 
\[ \dim(N_j \ominus N_{j-1}) \ge \dim(N_{j+1} \ominus N_j) \qfor 1 \le j < q .\]
Then $\rtsr(\T(\N)) = \infty$.
\end{cor}

\begin{cor}\label{C:omega2}
If a nest $\N$ has ordinal type at least $\omega^2$, then 
\[ \ltsr(\T(\N)) = \rtsr(\T(\N)) = \infty .\]
\end{cor}

\begin{proof}
By \cite[Theorem~2.2]{DLMR}, we have $\ltsr(\T(\N)) = \infty$.
By hypothesis, there is an increasing sequence $N_0<N_\omega < N_{2\omega} < \dots$
so that the intervals $N_{k\omega} \ominus N_{(k-1)\omega}$ are
all infinite dimensional for $k \ge 1$.
So by Corollary~\ref{C:decreasing}, we have $\rtsr(\T(\N)) = \infty$.
\end{proof}

\begin{rem}
 This takes care of most ordinals.
If $\N$ has infinite ordinal type less than $\omega^2$, then it can be 
written as a finite upper triangular matrix with the diagonal entries being nest algebras
of order type $\omega$ or finite, and upper triangular entries arbitrary.
A repeated application of Lemma~\ref{L:2x2} shows that $\rtsr(\T(\N))$
is just the maximum of the right stable ranks of the diagonal entries.
So the general problem is now reduced to studying nests of order type $\omega$
with finite rank atoms.
\end{rem}

\section{The Monotone Case and Beyond}

The situation simplifies considerably if the dimensions of the atoms 
is a monotone increasing sequence.  We will calculate the topological stable rank
in this case, and then extend the result to a chain of subnests, each
relatively monotone in the next.

\begin{thm}\label{T:monotone}
Let $\N=\{N_k, \H : k \ge0\}$ be a nest of order type $\omega$
with finite dimensional atoms of rank $n_k = \dim(N_k \ominus N_{k-1})$
for $k \ge 1$.
Suppose that the sequence $(n_k)$ is monotone increasing: \vspace{.3ex}
$n_{k+1} \ge n_k$ for $k\ge1$.
Define $d_j = \big| \{ k : 2^{j-1} < n_k \le 2^j \} \big|$ for $j \ge 0$.
\begin{enumerate}
\item  If $\sup d_j = \infty$, then $\rtsr(\T(\N)) = \infty$.
\item  If $\sup d_j < \infty$, then $\rtsr(\T(\N)) = 2$.
\end{enumerate}
\end{thm} 

\begin{proof}
Define $p_j = \sup\{ i : n_i \le 2^j\}$ for $j \ge 0$.  
Then $d_0=p_0$ and $d_j = p_j-p_{j-1}$ for $j \ge 1$.

If (i) holds, then given $q$, select some $j$ so that $d_j \ge 2^q$.
Let $k_0=p_{j-1}$ and 
\[ k_i = p_{j-1}+ 2^{q-1} + \dots + 2^{q-i} = p_{j-1} + 2^q - 2^{q-i} \qfor 1 \le i \le q .\]
Consider the chain $N_{k_0} < N_{k_1} < \dots < N_{k_q}$.
Then  $\dim(N_{k_i} - N_{k_{i-1}})$ lies in the range 
$(2^{q-i} 2^{j-1}, 2^{q-i} 2^j] = (2^{q+j-i-1}, 2^{q+j-i}]$.
In particular, these dimensions are monotone decreasing.
Therefore by Corollary~\ref{C:decreasing}, $\rtsr(\T(\N)) = \infty$.

If (ii) holds, then there is an integer $D$ so that $d_j \le D$ for all \mbox{$j\ge0$.}
For any $k \in (p_{j-1},p_j]$, we have $2n_k < 2^{j+1} < n_{k+2D}$. 
Therefore, \mbox{\cite[Theorem~2.11]{DLMR}} applies.
Hence $\rtsr(\T(\N)) = 2$.
\end{proof}

This applies to subexponential growth of the atom sizes:

\begin{cor} \label{C:subexp}
 If $(n_k)$ is monotone and $\liminf_{k\to\infty} n_k^{1/k} = 1$, then $\rtsr(\T(\N)) = \infty$.
\end{cor}

\begin{proof}
Fix $s$ and find $k>s^2$ so that $n_k^{1/k} < 2^{1/s}$.
Let $N = \lceil \frac{k}{s}\rceil$.  Then $p_N \ge k$.  Hence
\[ \sup d_j \ge \frac{p_N}N \ge \frac{k}{1+\frac{k}{s}} = \frac{s}{1+\frac s k} >  \frac{s}{1+\frac 1 s} > s-1.\]
As $s$ was arbitrary, this falls under case (i) of Theorem~\ref{T:monotone},
and so $\rtsr(\T(\N)) = \infty$.
\end{proof}

\begin{eg}
This answers Question 2 of \cite{DLMR}.  If $\N$ is a nest of order type $\omega$
with atoms of dimension $n_k = k$, then $\rtsr(\T(\N)) = \infty$.
This also provides a quite different proof that $\rtsr(\T(\N)) = \infty$ 
when $n_k = 1$ for all $k\ge1$.
\end{eg}

The following example motivates the subsequent analysis.

\begin{eg}
Consider the nest $\N$ of order type $\omega$ with atom dimensions given by the
sequence $(1), (1,2), (1,2,4), (1,2,4,8), \dots$.  This is not monotone.
However, we have grouped the atoms into segments on which they are monotone.
Consider the subnest $\M$ obtained by combining the groups in parentheses into single atoms,
so they have dimensions $1,3,7,15,\dots$.
This falls under case (ii) of Theorem~\ref{T:monotone}; so $\rtsr(\T(\M))=2$.
Let $\fJ = \ker\Delta$ be the kernel of the diagonal expectation onto $\fD(\M)$
as before.  Then $\rtsr(\tilde\fJ) = 2$ by Theorem~\ref{T:split}.

Consider $\fJ$ as an ideal of $\T(\N)$.  The quotient $\T(\N)/\fJ$ is isomorphic to
$\fB = \Delta(\T(\N)) = \bigoplus_{j=1}^\infty \T(\N_j)$ where $\N_j$ is the finite
nest with atoms of dimension $1,2,4,\dots,2^{j-1}$.
Each of these nest algebras is a quotient of the nest algebra $\T(\P)$
which has atoms of size $1,2,4,\dots$.
Now Theorem~\ref{T:monotone} shows that $\rtsr(\T(\P))=2$.
By Propositions~\ref{P:quotients} and \ref{P:direct sum}, it follows that
$\rtsr(\fB) = 2$.
Then by Theorem~\ref{T:split}, we see that $\rtsr(\T(\N)) = 2$.

Note that $\lim_{k\to\infty} n_k^{1/k} = 1$.  
So Corollary~\ref{C:subexp} is not valid for sequences which are not monotone.
\end{eg}

Since we are considering quotients obtained by expectations onto the diagonal
of a subnest of a nest, it will be useful to have such a variant of 
Theorems~\ref{T:monotone}.
In that proof, we used the quantitative lower bound estimate for part (i),
but could get away in part (ii) with the qualitative version.
Here we will need the estimates from Theorem~\ref{T:quantitative 2.11}.

{\samepage 
\begin{cor} \label{C:sum monotone}
Let $\N_i $ be nests for $i\ge1$ which are finite or order type $\omega$
and have only finite dimensional atoms of rank $n_{ik}$ for $1 \le k < K_i$,
where $K_i \in \bN \cup \{\infty\}$ is the cardinality of $\N_i$.
Let $\fB = \bigoplus_{i\ge1} \T(\N_i)$.
Suppose that the sequences $n_{ij}$ are monotone increasing in $j$ 
for each $i\ge1$.
For each $i\ge1$, define $d_{ij} = \big| \{k: 2^{j-1} < n_{ik} \le 2^j\} \big|$ for $j \ge 0$.
If $\sup_{i\ge1} K_i < \infty$, then $\tsr(\fB) = 1$.  Otherwise
\begin{enumerate}
 \item  if $\sup d_{ij} = \infty$, then $\rtsr(\fB) = \infty$.
 \item  if $\sup d_{ij} < \infty$, then $\rtsr(\fB) = 2$.
\end{enumerate}
\end{cor}
}

\begin{proof}
If $\sup_{i\ge1} K_i < \infty$, then $\tsr(\fB) = 1$ by Remark~\ref{R:finite nests}.
Otherwise the fact that the invertibles are not dense in $\T(\N)$ for any infinite
nest means that there can be no uniform control on $\rho_1(\T(\N_i),t)$ as 
the length of $\N_i$ tends to infinity.  Thus $\fB$ has $\rtsr(\fB)\ge2$.

If  $D = \sup d_{ij} = \infty$, then by Theorem~\ref{T:lower bound},
we obtain that 
\[ \sup_{i\ge1} \rho_n(\fA,t) = \infty \qforal n \ge1 \AND 0 < t \le \tfrac14 .\]
Hence $\rtsr(\fB) = \infty$.

Suppose that $D < \infty$. As in the proof of Theorem~\ref{T:monotone}, 
$J=2D$ satisfies the hypotheses of Theorem~\ref{T:quantitative 2.11} with $r=2$
for each nest $\N_i$.
We choose $p$ so that $2^p>5pJ$.  Then each of the nests $\N_i$ satisfy
$\rho_2(\T(\N_i),t) \le C t^{-5pJ-4}$ for a common constant $C$.
Hence by Proposition~\ref{P:direct sum}, $\rho_2(\fB,t) \le C t^{-5pJ-4}$
for $0<t\le1$.  Hence $\rtsr(\fB) = 2$.
\end{proof}

\begin{rem}
Suppose that $\N$ has a finite subnest $\M$ such that the 
atoms of $\N_i = \N \cap E_i$ have monotone increasing dimensions. 
Then the final atom of $\M$ has infinite rank; and so the dimensions of the 
atoms of $\N$ are monotone increasing for $k \ge k_0$ for some $k_0$.
By Lemma~\ref{L:2x2}, $\rtsr(\T(\N))$ depends only on this terminal sequence.
This is computed using Theorem~\ref{T:monotone}.
\end{rem}

\begin{defn}
If $\N = \{N_k,\H : k \ge 0 \}$ is a nest of order type $\omega$ 
with finite rank atoms of dimension $(n_k)$
and $\M = \{ N_{k_i},\H : i \ge 0\}$ is a subnest with $0=k_0 < k_1 < \dots$,
say that $\N$ is \textit{monotone relative to $\M$} if for each atom $E_i$ of $\M$,
the nest $E_i \cap \N$ has atoms of monotone increasing dimension.
The \textit{minimal relatively monotone subnest} of $\N$ is the smallest subnest
with this property, namely
\[ n_1 \le \dots \le n_{k_1} > n_{k_1+1} \le \dots \le n_{k_2} > n_{k_2+1} \le \dots .\]
\end{defn}

\begin{thm} \label{T:k-monotone}
Let $\N$ be a nest of order type $\omega$ with finite rank atoms.
Suppose that there is a sequence of subnests 
$\N = \N_0 \supset \N_1 \supset \dots \supset \N_k$
such that for each $j \ge 1$, either
\begin{enumerate}
\item[(1a)] $\N_j$ is of finite index in $\N_{j-1}$, or
\item[(1b)] $\N_{j-1}$ is monotone relative to $\N_j$; and
\item[(2)] the atom dimensions of $\N_k$ are monotone increasing.
\end{enumerate}
Then $\rtsr(\T(\N)) \in \{ 2,\infty\}$; and there is an algorithm to compute it.
\end{thm}

\begin{proof}
By Theorem~\ref{T:chain of nests}, 
\[
 \rtsr(\T(\N)) = \max\big\{ 
 \rtsr(\T(\N_k)), \rtsr(\Delta_{\N_j}(\T(\N_{j-1}))), 
 1 \le j \le k \big\} .
\]
First apply Theorem~$\ref{T:monotone}$ to compute $\T(\N_k)$.
Then apply Corollary~\ref{C:sum monotone} to each of the algebras 
$\fB_{j-1} = \Delta_{\N_j}(\T(\N_{j-1}))$. 
All of these values lie in $\{2,\infty\}$.
So the result follows.  Moreover, this provides an algorithm for the calculation.
\end{proof}

\begin{rem}
In the proof  of part (i) of Theorem~\ref{T:monotone}, 
the monotonicity is not the critical issue.
What was used is that there are arbitrarily long intervals of $\N$ on which the
dimensions of the atoms are all comparable within a factor of 2.
Clearly, 2 can be replaced by any constant $\gamma>1$ and a suitable modification
of the argument still shows that $\rtsr(\T(\N))=\infty$.

Call a sequence $(n_k)$ of positive integers \textit{almost increasing}
if there is a constant $\gamma \ge 1$ so that $n_j \le \gamma n_k$ for all $1 \le j < k$.

One can prove variations of Theorems~\ref{T:monotone} and \ref{T:k-monotone}
replacing monotone by almost monotone.  
This is also the case for Corollary~\ref{C:sum monotone} provided that
one assumes a common constant $\gamma$ for all summands.
\end{rem}

\section{Topological stable rank of nest algebras}\label{S:nests}

In this section, we characterize the right (left) topological stable rank
of an arbitrary nest algebra.  The results of the previous section will provide
a method for computing it in a more practical sense.

The following definition captures the invariant used in Corollary~\ref{C:decreasing}.

\begin{defn}
If $\N$ is a nest, let $\beta(\N)$ be the supremum of integers $q$ for which
there is a chain $N_0<N_1<\dots<N_q$ in $\N$ with
\[ \dim(N_j\ominus N_{j-1}) \ge \dim(N_{j+1}\ominus N_j) \qfor  1 \le j \le q .\]
\end{defn}

\begin{lem} \label{L:beta1}
If $\N$ is a nest with $\beta(\N) < \infty$, then $\N$ has ordinal type $\alpha < \omega^2$
and only finitely many infinite rank atoms.
\end{lem}

\begin{proof}
Suppose that $\N$ has an infinite decreasing sequence $(N_k)_{k \ge 0}$,
$N_k > N_{k+1}$ for all $k \ge 0$.  Let $E_k = N_{k-1} \ominus N_k$ for $k \ge1$.
There are two case to consider.  
If infinitely many of these intervals have infinite dimension,  then by dropping to a
subsequence, we may suppose that all the intervals have infinite dimension.
Hence there are arbitrarily long sequences with the intervals all having the same infinite
dimension.  So $\beta(\N) = \infty$.
Otherwise, only finitely many of these intervals are infinite dimensional.
By deleting the initial terms of our sequence, we may then suppose
that all of the intervals have finite dimension.
It is routine to drop to a subsequence $(N_{k_i})_{i\ge0}$ so that
the intervals \mbox{$F_i = N_{k_{i-1}} \ominus N_{k_i}$} have monotone increasing dimension.
Thus the sequence $N_{k_q} < N_{k_{q-1}} < \dots < N_{k_0}$ is a chain of decreasing
dimension for arbitrary $q\ge 1$.  So again $\beta(\N) = \infty$.

Arguing as in \cite[Theorem~2.4]{DLMR}, we see that any nest which is not order
isomorphic to an ordinal will have an infinite decreasing sequence.
Also, when $\N$ is order isomorphic to an ordinal $\alpha \ge \omega^2$,
we argue as in Corollary~\ref{C:omega2}  with the sequence 
$N_0 < N_\omega < N_{2\omega} < \dots$
to show that $\beta(\N) = \infty$.
Finally, if $\N$ has ordinal type $\alpha < \omega^2$ but has infinitely many infinite
rank atoms, then it is also easy to construct an infinite increasing sequence of
infinite intervals; whence $\beta(\N) = \infty$.
\end{proof}

By Corollary~\ref{C:decreasing}, $\beta(\N) = \infty$ implies that $\rtsr(\T(\N)) = \infty$;
as does the existence of any infinite dimensional atom.
We will establish the converse, and at the same time we will
show that $\rtsr(\T(\N)) = 2$ is the only other possibility.

{\samepage
\begin{lem} \label{L:beta2}
Let $\N=\{N_k,\H : k \ge 0\}$ be a nest of order type $\omega$ with 
atoms of finite rank $n_k$ for $k \ge 1$.
Suppose that $\beta(\N) = q < \infty$. 
Let $\N'$ be the minimal relatively monotone subnest of $\N$.
\begin{enumerate}
\item If $q=1$, then $\N' = \{0,\H\}$ and $n_k > \sum_{i=1}^{k-1} n_i$ for all $k \ge 2$.
\item If $q>1$, then $\beta(\N') \le q-1$.
\end{enumerate}
\end{lem}
}
\begin{proof}
Let $A_k = N_k \ominus N_{k-1}$ be the atoms of $\N$.
First suppose that $\beta(\N) = 1$.
Consider the chain $N_0 < N_{k-1} < N_k$ for $k \ge 2$.
Since this cannot have decreasing dimension of intervals,
it follows that 
\[ \sum_{i=1}^{k-1} n_i = \dim N_{k-1} < \dim A_k = n_k .\]
In particular, $(n_k)$ is monotone increasing, so that $\N' = \{0,\H\}$.

Suppose that $q \ge 2$.  Let $\N' = \{ N_{k_i},\H : i \in \I\}$
(which may be finite).
Suppose that $\N'$ has a chain $N'_0 < N'_1 < \dots < N'_p$ so that
$E_j = N'_j \ominus N'_{j-1}$ have decreasing dimension.
These dimensions are all finite, so $N'_p \ne \H$.
The interval $E_p$ dominates the atom $A' = N'_p \ominus N^{\prime -}_p$.
Also $A'\cap \N$ is a maximal interval of $\N$ on which the 
sequence $n_i$ is monotone increasing.  This means that
the next atom of $\N$, say $A_{k_0}$, must have smaller
rank.  In particular, $\dim A_{k_0} < \dim A' \le \dim E_q$.
So the sequence $N'_0 < N'_1 < \dots < N'_p < N_{k_0}$ in $\N$
would be a sequence of decreasing dimension of length $p+1$.
Thus $p+1 \le q$; so that $\beta(\N') \le q-1$.
\end{proof}

The following result completely characterizes the right topological stable rank of
nest algebras.  In particular, it answers Question 1 of \cite{DLMR} by showing that
the only possible values are $2$ and $\infty$.

{\samepage 
\begin{thm} \label{T:nests}
For a nest $\N$ on separable Hilbert space, the following are equivalent:
\begin{enumerate}
\item $\rtsr(\T(\N)) = 2$.
\item $\rtsr(\T(\N)) < \infty$.
\item $\beta(\N) < \infty$ and $\N$ has no infinite rank atoms.
\end{enumerate}
\end{thm}
}

\begin{proof}
Clearly (i) implies (ii).  
Corollary~\ref{C:decreasing} shows that (ii) implies (iii).
Suppose that $\beta(\N) < \infty$ and that all atoms are finite dimensional.
By Lemma~\ref{L:beta1}, this shows that $\N$ is of ordinal type $\alpha < \omega^2$.
So $\alpha \le n\omega$ for some integer $n$.
That means that $\T(\N)$ is a finite upper triangular matrix
with diagonal entries $\T(\N_i)$ for nests of order type $\omega$ or finite.
Clearly $\beta(\N_i) \le \beta(\N) < \infty$.
By Lemma~\ref{L:2x2}, it suffices to show that $\rtsr(\T(\N_i)) = 2$.
Thus the problem is reduced to nests of order type $\omega$ with finite
rank atoms and finite $\beta$.

There is a sequence of subnests 
$\N = \N_0 \supset \N_1 \supset \dots \supset \N_p$ where
each $\N_{j+1}$ is the minimal relatively monotone subnest of $\N_j$.
By Lemma~\ref{L:beta2}, $\beta(\N_{j+1}) < \beta(\N_j)$.
So we may proceed until $\N_p$ is a finite nest.
Therefore $\N_{p-1}$ has monotone increasing atom dimensions
at least from some point on.

Therefore we satisfy the hypotheses of Theorem~\ref{T:k-monotone}.
Let $d_j$ be defined for the nest $\N_{p-1}$ as in Theorem~\ref{T:monotone}.
The proof of Theorem~\ref{T:monotone} is accomplished by showing
that $\sup d_j = \infty$ implies that $\beta(\N) = \infty$.
Hence we conclude $\sup d_j < \infty$, and so $\rtsr(\T(\N_{p-1})) = 2$.
Therefore we have
\[ \rtsr(\T(\N)) = \max\{2, \rtsr(\Delta_{\N_j}(\T(\N_{j-1}))) : 1 \le j \le p-1\} , \]
where $\Delta_{\N_j}$ is the expectation onto the diagonal of $\N_j$.
When $E$ is an atom of $\N_j$,  $E \cap \N_{j-1}$ consists of atoms of
increasing dimension. 
So Corollary~\ref{C:sum monotone} applies.
Again there are two possibilities, and the first, $\sup d_{ij} = \infty$
leads to the conclusion that $\beta(\N_{j-1}) = \infty$, contrary to fact.
Thus $\rtsr(\Delta_{\N_j}(\T(\N_{j-1}))) = 2$ for every $1 \le j < p-1$.
Therefore $\rtsr(\T(\N)) = 2$.
\end{proof}

{\samepage 
\begin{cor} \label{C:2 or infinity}
For any nest $\N$ on an infinite dimensional Hilbert space, 
\[ \rtsr(\T(\N)) \in \{2,\infty\} \qand  \ltsr(\T(\N)) \in \{2,\infty\} \]
and 
\[ \max\{\rtsr(\T(\N)), \ltsr(\T(\N)) \} = \infty .\]
\end{cor}
}

\begin{proof}
Theorem~\ref{T:nests} shows that $\rtsr(\T(\N))$ is either $2$ or $\infty$.  
The same follows for left topological stable rank
because $\ltsr(\T(\N)) = \rtsr(\T(\N^\perp))$.  
The last identity 
is established in \cite[Corollary~2.5]{DLMR}.
\end{proof}

Given a specific nest of order type $\omega$,  it may still require some
non-trivial effort to calculate $\rtsr(\T(\N))$.  The proof explains the algorithm.
One successively constructs the minimal monotone subnests 
\[ \N \supset \N_1 \supset \N_2 \supset \cdots .\]
If $\rtsr(\T(\N)) = 2$, this process will terminate.
Thus if this sequence never reaches the trivial nest, one must have $\rtsr(\T(\N)) = \infty$.
On the other hand, it this is a finite chain, then Theorem~\ref{T:k-monotone}
applies, and the steps in the chain are analyzed by the use of 
Theorem~\ref{T:monotone} and Corollary~\ref{C:sum monotone}.

\section{Partial Matrix algebras}\label{S:matrix}

Rieffel \cite[Theorem~6.1]{R} establishes a very nice, precise result about 
full matrix algebras over a Banach algebra $\fA$, namely
\[  \rtsr(\fM_n(\fA))= \Big\lceil \frac{\rtsr(\fA)-1}{n} \Big\rceil + 1 .\] 
We are interested in comparing two nest algebras with similar growth
which are not exactly related by multiplicity.  So we need a modification
of Rieffel's argument for a larger class of matrix-like algebras over $\fA$
with a corresponding weakening of the bounds on the stable rank.

\begin{defn}
Say that a set of \textit{$n \times n$ partial matrix units} is a collection of idempotents
$P_i$, $1\le i \le n$, such that $\sum_{i=1}^n P_i = I$ 
and operators $U_i, V_i$ for $2 \le i \le n$ satisfying
\[ U_i = P_i U_i P_1, \qquad V_i = P_1V_iP_i \qand U_iV_i=P_i \qfor 2 \le i \le n. \]
\end{defn}

The reason these are partial matrix units is that we do not require that $V_iU_i = P_1$.
Now $V_iU_i = E_i$ will always be an idempotent such that $E_i = E_iP_1 = P_1E_i$.
Let $E_1 = P_1$.  For notational convenience, we also set $V_1=U_1=P_1$.
Set $\fA_1$ to be the Banach algebra $P_1\fA P_1$ with unit $P_1$.
Then we can consider $\fA$ as a \textit{partial matrix algebra over $\fA_1$}.
Indeed, it is easy to see that $\fA$ is isomorphic to the Banach algebra 
$E\,\fM_n(\fA_1) E$, where $E = \sum_{i=1}^n \oplus E_i$.

For convenience, we will normalize the partial matrix units so that $\|V_i\|=1$.
Identify $\fA_1$ with $P_1\fA P_1$, which we think of as the $1,1$ entry of $\fA$.
Then for $A\in\fA$, 
\[ P_i A P_j = U_iV_i A U_jV_j = U_i (P_1 V_i A U_j P_1) V_j .\]
Set $a_{ij} = P_1 V_i A U_j P_1$ considered as an element of $\fA_1$.
Then $A$ may be thought of as an $n\times n$ matrix with coefficients
$A_{ij} = U_i a_{ij} V_j$.

The main result of this section is the following result.

\begin{thm} \label{T:partial matrix units}
Let $\fA$ be a Banach algebra with partial matrix units, and let $\fA_1 = P_1\fA P_1$.
Then
$\rtsr(\fA_1) \le (\rtsr(\fA) - 1)n + 1$ and hence
\[
  \left\lceil \frac{\rtsr(\fA_1) - 1}{n} \right\rceil + 1 \le \rtsr(\fA).
\] 
Suppose that $E_i = P_1$ for $1 \le i \le m$.  
If $\fA_1$ is completely finite, then
\[
 \rtsr(\fA) \le \left\lceil \frac{\rtsr(\fA_1)-1}{m} \right\rceil + 1 . 
\]
\end{thm}

\begin{proof}
Let $\rtsr(\fA) = p$ and $L = \max\{\|U_i\|\,\|V_i\| : 2 \le i \le n \}$.
We first establish that $\rtsr(\fA_1) \le (p-1)n+1$ in the same  \vspace{.5ex}
manner as in Rieffel's proof.
This will immediately imply that 
$\big\lceil \frac{\rtsr(\fA_1)-1}{n} \big\rceil + 1 \le \rtsr(\fA)$. \vspace{.5ex}

Start with $a = \big[ a_1\ \dots\ a_{(p-1)n+1} \big] \in \R_{(p-1)n+1}(\fA_1)$.
Define $p$ elements of $\fA$ by
\[
 B_s = \sum_{j=1}^n a_{(s-1)n+j} V_j = 
 \begin{bmatrix} a_{(s-1)n+1}V_1 & \dots & a_{sn}V_n\\ 0 & \dots & 0\\
 \vdots & \ddots & \vdots\\ 0 & \dots & 0 \end{bmatrix}
 \qfor 1 \le s < p
\]
and
\[
 B_p = a_{(p-1)n+1} + \sum_{j=2}^n P_j =
 \begin{bmatrix} a_{(p-1)n+1} & 0&\dots&0\\
 0& P_2 & \dots &0\\
 \vdots&\vdots&\ddots&\vdots\\
 0&0&\dots&P_n \end{bmatrix} .\qquad\quad\strut
\]
So $B = \big[ B_1 \ \dots \ B_p \big] \in \R_p(\fA)$.
We may think of $B$ as having the form 
$B = \begin{bmatrix}aV & 0\\0 & I_{n-1}\end{bmatrix}$,
where $V = \diag\big( V_{j\!\!\mod n} : 1 \le j \le (p-1)n+1 \big)$ and
$I_{n-1} = \sum_{j=2}^n P_j$ is the identity on the 
last $n-1$ blocks of the matrix.

Given $0 < t \le 1$, there is a perturbation $B'$ of $B$ of norm at most $\frac t{L+1}$
which is right invertible; say 
\[ B' = \begin{bmatrix}a' & b\\c & I_{n-1}+d\end{bmatrix} \]
with 
\[
 \|B'-B\| =
 \left\| \begin{bmatrix}a' - aV & b\\c & d\end{bmatrix} \right\| 
 \le \frac t{L+1}
\]
and  right inverse $C$. 
Define an invertible matrix
\[
 D = \begin{bmatrix} P_1 & b\\ 0 & I_{n-1}+d \end{bmatrix} 
 \quad\text{with}\quad
 D^{-1} = \begin{bmatrix} P_1 & -b(I_{n-1}+d)^{-1}\\ 0 & (I_{n-1}+d)^{-1}\end{bmatrix} .
\]
Now 
\[
 D^{-1}B' = \begin{bmatrix}a' - b (I_{n-1}+d)^{-1}c & 0\\
 (I_{n-1}+d)^{-1}c & I_{n-1} \end{bmatrix} .
\]
This is right invertible by $CD$.

The row $a'' := (a'-aV) - b (I_{n-1}+d)^{-1}c$ has norm at most
\[
 \|a'-aV\| + \frac{\|b\|\,\|c\|}{1-\|d\|}  \le 
 \frac t {L+1} +  \frac {t^2} {(L+1)^2} \frac {1}{1-\frac t{L+1}} 
 \le \frac t L .
\]
Set $U = \diag \big( U_{j\mod n} : 1 \le j \le (p-1)n+1 \big)$. 
Then 
\[ a'' = (a''U)V = yV .\]
where $y = a''U$ belongs to $\R_n(\fA_1)$.
So $\|y\| \le \|a''\|\,\|U\| \le t$.
Thus we have 
\[ a' - b (I_{n-1}+d)^{-1}c = (a+y)V .\]

Let $X$ denote the column consisting of the first \mbox{$(p\!-\!1)n+1$} entries
of the first column of $CD$. 
We may write this as $X = \big[ U_{j\!\!\mod n}\, x_j \big]^t$ for $x_j \in \fA_1$.
Then 
\[ P_1 = \big( a' - b (I_{n-1}+d)^{-1}c \big) X .\]
Note that $VX = \big[ E_{j\!\!\mod n}\, x_j \big]^t$ belongs to 
$\C_{(p\!-\!1)n+1}(\fA_1)$.
Finally observe that 
\[ (a+y)(VX) = ( aV+a'' )X = P_1 .\]

\bigbreak
We now turn to the second inequality.
If $E_i=P_1$ for $1 \le i \le m$, then the corner $m\times m$ 
submatrix of $\fA$ is isomorphic to $\fM_m(\fA_1)$, and this algebra \vspace{.3ex}
has $\rtsr(\fM_m(\fA_1))= \Big\lceil \frac{\rtsr(\fA_1)-1}{m} \Big\rceil + 1$
by Rieffel's Theorem.
So we may consider $\fA$ as a partial matrix algebra over 
$\fM_m(\fA_1)$ instead.
Since $\fA_1$ is completely finite, so is $\fM_m(\fA_1)$ 
by Corollary~\ref{matrix comp finite}.
This reduces the problem to the case of $m=1$.

Let $\rtsr(\fA_1) = q$. Fix $0 < t \le 1$.
Consider a row $R = \big[ R_1 \ \dots \ R_q \big]\in \R_q(\fA)$ with $\|R\| \le 1$.
Write $R_i = \begin{bmatrix} a_i & b_i\\ c_i & d_i \end{bmatrix}$ with respect to the
decomposition $P_1\H \oplus \sum_{i=2}^n P_i \H$.
By doing the `canonical shuffle', we may rearrange this as
$R =  \begin{bmatrix} A & B\\ C & D \end{bmatrix}$ where
each entry is a $1\times q$ row; 
e.g. $A = \big[ a_1 \ \dots \ a_q \big] \in \R_q(\fA_1)$.
This has a $t/n$ perturbation $A'$ with right inverse $X \in \C_q(\fA_1)$.

Since $\fA_1$ is completely finite, there is an invertible matrix $W_1$ in $\fM_q(\fA_1)$
with first row $A'$ such that the first column of $W_1^{-1}$ is $X$.
Let $R' =  \begin{bmatrix} A' & B\\ C & D \end{bmatrix}$,  \vspace{.5ex}
and let $\tilde B$ be the $q\times q$ matrix with first row equal to $B$ 
and the other rows zero.  Compute
\begin{align*}
 \ol{R} &= 
 \begin{bmatrix} A' & B\\ C & D \end{bmatrix}\, 
 \begin{bmatrix} W_1^{-1} & -W_1^{-1} \tilde B\\ 0 & I_{n-1}^{(q)} \end{bmatrix}
 = \begin{bmatrix} J_1 & 0\\ C_1 & D_1 \end{bmatrix}
\end{align*}
where $J_i = \big[ P_i\ 0\ \dots\ 0 \big]$. 
Since $R'$ has been multiplied by an invertible matrix to get $\ol{R}$,
it is close to a right invertible element if and only if $ \ol{R}$ is.
Indeed, if 
\[
 w_1 = \norm{ \begin{bmatrix} W_1 & \tilde B\\ 0 & I_{n-1}^{(q)} \end{bmatrix} } \le \|W_1\|+1 ,
\]
then a $t/nw_1$ perturbation of $\ol{R}$ translates back to a $t/n$ perturbation of $R$.

The plan now is to use the matrix units to move each diagonal entry
into the $1,1$ entry one at a time and apply the same procedure as 
above to the result each time.  We illustrate this for $n=3$, and the
reader will see how to set it up as an induction.

First we deal with the $P_2$ block.  
Let us write 
\[
 C_1 =  \begin{bmatrix} C_{21}\\C_{31} \end{bmatrix} 
 \qand  
 D_1 =  \begin{bmatrix} D_{22}&D_{23}\\D_{32}&D_{33} \end{bmatrix}
\]
where $C_{i1} \in \R_q(P_i \fA P_1)$ and $D_{ij} \in \R_q(P_i \fA P_j)$.
Compute
\begin{align*}
 &\begin{bmatrix}
  P_1 \!-\! E_2 & V_2 & 0 \\
  U_2 &  0 & 0 \\ 
  0 & 0 & P_3 
 \end{bmatrix} 
 \begin{bmatrix}
  J_1 & 0 & 0\\
  C_{21} & D_{22} & D_{23} \\ 
  C_{31} & D_{32} & D_{33}
 \end{bmatrix}
 \begin{bmatrix}
  (P_1 \!-\! E_2)^{(q)} & V_2^{(q)} & 0 \\ 
  U_2^{(q)} &  0 & 0 \\ 
  0 & 0 & P_3^{(q)} 
 \end{bmatrix}  \\
 &\ = 
 \begin{bmatrix}
  (P_1 \!-\! E_2)J_1 + V_2C_{21}(P_1\!-\!E_2)^{(q)} + V_2 D_{22} U_2^{(q)}  
  & V_2 C_{21} V_2^{(q)} & V_2 D_{23} \\ 
  0 &  J_2  & 0 \\ 
  C_{31}(P_1\!-\!E_2)^{(q)} + D_{32} U_2^{(q)}  & C_{31}V_2^{(q)} & D_{33}
 \end{bmatrix}  \\ &\ = 
 \begin{bmatrix}
  A_2 & B_{12} & B_{13} \\ 
  0 &  J_2  & 0 \\ 
  C'_{31}  & C'_{32} & D_{33}
 \end{bmatrix} 
\end{align*}
Multiplying on either side by an invertible matrix, such as here, 
takes the set of right invertible rows onto itself.  
So it suffices to find a small perturbation of this new matrix
which is right invertible.
Note that the matrices used above are symmetries ($S^2=I$) 
and their norms are bounded by 
$\max\{\|U_2\|,\|V_2\|\} + \|P_1 - V_2U_2\| \le 2L+1$.
So a $t(nw_1)^{-1}(2L+1)^{-2}$ perturbation here translates
back to a $t/n$ perturbation of the original $R$.

We can now find a small perturbation of the $1,1$ entry, 
say $A'_2$, with right inverse $X_2$.
Let $W_2$ be an invertible matrix in $\fM_q(\fA_1)$ with first row $A_2$
and such that $W^{-1}$ has first row $X_2$.
Let $\tilde B_{1j}$ denote the $q\times q$ matrix with first row $B_{1j}$
and the remaining rows equal to $0$.
Again we multiply
\[
 \begin{bmatrix}
  A'_2 & B_{12} & B_{13} \\ 
  0 &  J_2  & 0 \\ 
  C'_{31}  & C'_{32} & D_{33}
 \end{bmatrix} 
 \begin{bmatrix}
  W_2^{-1} & -W_2^{-1} \tilde B_{12} & -W_2^{-1} \tilde B_{13} \\ 
  0 &  P_2^{(q)}  & 0 \\ 
  0 & 0 & P_3^{(q)}
 \end{bmatrix} 
 =
 \begin{bmatrix}
  J_1 & 0 & 0 \\ 
  0 &  J_2  & 0 \\ 
  C''_{31}  & C''_{32} & D_{33}
 \end{bmatrix} 
\]
This has a small perturbation which is right invertible only if the
previous matrix had such a perturbation.

Clearly this procedure may be repeated until one obtains
the matrix $\diag ( J_1, J_2, \dots, J_n)$. After undoing the canonical shuffle,
this is just the row $\big[ I\ 0 \ \dots \ 0 \big]$, which is evidently right invertible.
\end{proof}

\begin{rem}
There is a quantitative version of this theorem.
For the first part, we can obtain
\[
 \rho_{(p-1)n+1}(\fA_1,t) \le 2 \rho_p(\fA,\tfrac t{L+1}) .
\]
To obtain a quantitative version of the second part, we need to require
that $\fA_1$ is uniformly completely finite.
The details are left for the interested reader.
\end{rem}

This allows us to show that for nests, it is just the relative size of the atoms which 
is relevant in the determination of topological stable rank.

\begin{cor} \label{C:atom inequality}
Let $\N$ and $\M$ be nests of order $\omega$ with finite rank atoms of 
dimensions $n_k$ and $m_k$ for $k\ge1$.  
If there are positive constants $0 < c \le d$ such that $cn_k \le m_k \le dn_k$,
then 
\[ \rtsr(\T(\N)) = \rtsr(\T(\M)) .\]
\end{cor}

\begin{proof}
First assume that $c=1$ and $q=\lceil d \rceil$.
The atoms of $\M$ may be split into $q$ pieces, with the first having rank $n_k$
and the remainder of rank at most $n_k$.  Then it is evident that $\T(\M)$
is unitarily equivalent to a $q\times q$ partial matrix  algebra over $\T(\N)$
with $1,1$ entry equal to $\T(\N)$.
By Theorems~\ref{T:partial matrix units} and \ref{T:extend_nest},
\[ \left\lceil \frac{\rtsr(\T(\N)) - 1}q \right\rceil \le \rtsr(\T(\M)) \le \rtsr(\T(\N)) .\]
Since $\rtsr(\T(\N)) \in \{2,\infty\}$ by Corollary~\ref{C:2 or infinity},
it follows that the left hand side also equals $\rtsr(\T(\N))$;
and thus they are equal.

For the general case, let $\L$ be the nest of order type $\omega$
with atoms of rank $l_k = \min\{n_k,m_k\}$.
Then $l_k \le n_k, m_k \le q l_k$ for $q = \lceil \max\{c^{-1}, d\} \rceil$.
Hence 
\[  \rtsr(\T(\N)) = \rtsr(\T(\L)) = \rtsr(\T(\M)) . \qedhere \]
\end{proof}

\medskip\textit{Acknowledgements.} 
The second author would like to thank the University of Waterloo for its
hospitality during the period when this paper was written.



\end{document}